\documentclass[12pt]{article}
\usepackage[english]{babel}
\usepackage{epsfig}
\usepackage{amsfonts}
\usepackage{amssymb}
\usepackage{amsmath}
\usepackage{amsthm}
\usepackage{latexsym}
\usepackage{graphicx}
\usepackage{bm}
\usepackage{tabularx}
\usepackage{booktabs} 
\usepackage{enumerate}
\usepackage{caption}
\usepackage{wrapfig}
\usepackage{dsfont}

\usepackage{mathtools}  
\usepackage{eqnarray}
\usepackage{enumitem}

\usepackage[titletoc]{appendix}

\usepackage[numbers]{natbib}
\usepackage{comment}
\usepackage{array} 
\usepackage{bbm} 
\usepackage{url}
\usepackage{subfig}
\usepackage{hyperref}
\usepackage{color}
\usepackage{float}

\usepackage{hyperref}
\hypersetup{colorlinks=true, citecolor=blue, linkcolor=red, urlcolor=green}
 
\usepackage{bbm}
\usepackage{url}
\usepackage{color}
\usepackage{xcolor}
\usepackage[section]{placeins}

\theoremstyle{plain}
\newtheorem{theorem}{Theorem}[section]
\newtheorem{lemma}[theorem]{Lemma}

\newtheorem{proposition}[theorem]{Proposition}
\newtheorem{conjecture}[theorem]{Conjecture}

\theoremstyle{definition}

\theoremstyle{remark}
\newtheorem*{remark}{Remark}

\allowdisplaybreaks

\title{Overlap Times in Tandem Queues: Identically Distributed Station Case}
\author{ 
Ruici Gao
\\
School of Operations Research and Information Engineering
\\
Cornell University
\\
{rg585@cornell.edu} \\ \\
Jamol Pender \footnote{Corresponding Author}
\\
School of Operations Research and Information Engineering
\\
Cornell University
\\
{jjp274@cornell.edu}
}

\usepackage{natbib}
\usepackage{graphicx}
\usepackage{amsmath}

\begin{document}

\maketitle

\begin{abstract}
In this paper, we investigate overlap times in a two-dimensional infinite server tandem queue. Specifically, we analyze the amount of time that a pair of customers spend overlapping in any station of the two dimensional tandem network. We assume that both stations have independent and identically distributed exponential service times with the same rate parameter $\mu$. Our main contribution is the derivation of the joint tail distribution, the two marginal tail probabilities, the moments of the overlap times and the tail distribution of the sum of the overlap times in both stations. Our results shed light on how customers overlap downstream in serial queueing systems.  
\end{abstract}

\section{Introduction}


Tandem queues are a mathematical and analytical framework used to model and analyze complex queueing systems that consist of multiple sequential queues connected in series. In a tandem queue, customers or entities move from one queue to the next in a linear fashion, typically following a first-come, first-served (FCFS) discipline. This concept is particularly useful in various fields, including operations research, computer science, telecommunications, healthcare and manufacturing, to understand and optimize the performance of interconnected systems.

Applications of tandem queues are diverse and widespread. They are commonly used in telecommunications networks to model call routing and processing, where customers move through different stages of call handling, such as call setup, voice processing, and billing. They also model the transmission of data packets from one network node to another, where packets are processed in a sequential manner, see for example \citet{boxma1984m, ziedicnvs1993tandem, massey1993networks, dean1995spanning, daigle2005queueing, le2008tandem, pender2016analysis, palomo2021learning}. In manufacturing, tandem queues can help analyze production lines with multiple workstations, where products pass through various stages of assembly or testing, see for example \citet{avi1993servers, tsiotras1992wip}. Additionally, tandem queues find applications in computer systems, like web servers, where user requests are processed in a sequence of tasks, each representing a queue, and understanding the performance of such systems is crucial for efficient resource allocation and service quality optimization. 

In addition, tandem queues play an essential role in transportation systems, where passengers move from one vehicle to another. For example, in an airport, passengers may have to move from one terminal to another, and they may have to wait in queues to board a plane, see for example \citet{dunlay1978tandem}. In healthcare systems, tandem queues are used to model the movement of patients from one department to another. In hospitals, patients may have to wait in queues at each department, from the emergency department to the surgery department.  Overall, tandem queues are versatile and widely applicable, making them an essential modeling tool in many practical applications. Understanding and utilizing tandem queues can help optimize complex systems, reduce waiting times, and improve efficiency, making them an essential tool for managers and decision-makers in a variety of areas.

In this paper, we aim to combine the analysis of tandem queues with that of overlap times as seen in \citet{palomo2021overlap}.  Despite the significant practical relevance of overlap times in queueing networks, there is currently a gap in the literature regarding this topic. Existing research has mostly focused on one-dimensional queueing systems, as seen in studies such as \citet{palomooverlap, kang2021queueing, palomo2021overlap, ko2022overlapping, palomo2023overlap, xu2023queueing, ko2023number}. One exception to this is \citet{kang2021queueing}, which explores overlaps in a two dimensional priority queue.  To begin to address this gap in the literature, we present a novel study that explores overlap times in two-dimensional infinite server tandem queues. We focus on the two dimensional case since the two dimensional case is quite challenging and understanding higher dimensions will require new and more complicated techniques. Even in the two dimensional setting, we specifically analyze the overlap times in the Poisson arrival and exponential service setting, since it requires keeping track of overlaps in multiple queues simultaneously. 



\subsection{Application to Quantum Computing}

Despite the many applications of tandem queues and overlap times, we are partially inspired to analyze these systems since they arise in the area of quantum computation and quantum information processing.  One significant challenge in quantum computing is solving the quantum decoherence problem.  Quantum decoherence is a phenomenon that occurs in quantum systems when their delicate quantum properties, such as superposition and entanglement, are lost due to interactions with their environment. These interactions can include interactions with other particles, radiation, or any external factors that disturb the coherence of the quantum system. As a result of decoherence, a quantum system's behavior becomes more classical, thus making it behave more predictably according to classical physics.

Maintaining the coherence of qubits (quantum bits) is crucial for performing complex quantum computations and preserving the advantages of quantum information processing.  When two qubits are entangled or share an overlap in their quantum states, they become intertwined in a manner that makes them more resistant to certain types of decoherence. This entanglement or overlap can be harnessed to create logical qubits that are more robust against individual qubit errors. By exploiting the correlation between the states of the logical qubits, errors affecting one physical qubit can be detected and corrected using the information from other qubits in the logical system.  In essence, the overlap of two qubits, often achieved through entanglement, can prevent decoherence by creating redundancy and allowing for error correction. If one of the qubits experiences a decoherence-induced error, the information contained in the overlap or entanglement can be used to identify and correct the error, effectively restoring the original quantum state.

Overlap times also play a role in quantum interference, where the probability amplitudes of different quantum states can add constructively or destructively, leading to outcomes that are significantly different from classical probabilities. Quantum algorithms, like Shor's algorithm for factoring large numbers and Grover's algorithm for searching unsorted databases, leverage t his interference to achieve exponential speedup compared to their classical counterparts. 

Thus, it is critical in quantum systems to know the probability that qubits will overlap and entangle in order to speed up computations and realize quantum advantages.  Moreover, like bits, qubits are not isolated to just one quantum gate.  Qubits can travel through many gates to perform quantum computations.  The results in this work regarding tandem queues will allow quantum researchers to understand how qubits overlap in different quantum gates and understand the probability of having qubits entangle and realize quantum advantages.  

Fortunately, we are not the first to use concepts of queueing theory to analyze quantum systems, see for example \citet{dai2020quantumthesis, dai2020quantum,jagannathan2019qubits, mandalapu2022classical,  mandayam2019classical, mandayam2020classical,   varma2023erasure}.  However, most of this work analyzes the capacity of queue-channels and calculates the fraction of qubits lost because of decoherence and erasure.  Our work considers, the explicit interaction between qubits and measures their ability to entangle by understanding the probability of them overlapping a certain amount of time in a quantum circuit.  Moreover, we consider a network of queues, which is different from the previous literature as well.  
  
By viewing quantum gates as stations and qubits as arrivals, we can apply the results of this work to advance quantum computing. Analyzing the joint density of qubit overlap times reveals insights into circuit bottlenecks. Longer overlap times between qubit pairs lead to extended entanglement and coherent superposition, offering a strategy for fortifying circuits against decoherence.

\subsection{Contributions of the Paper}
\begin{itemize}
    \item We derive the joint and marginal distributions for the overlap time for customers that are exactly $k$ customers apart in a two dimensional infinite server tandem network queue with identical exponential service distributions. We also compute the tail distribution of the sum of the overlap times.
    \item We derive the joint and marginal distributions of the overlap time for different pairs of customers in both stations as well.  
    \item We compute the mean, variance and covariance for the overlap times in both stations. 
    \item We use simulation to verify our results.  
\end{itemize}

\subsection{Organization of Paper}
In Section \ref{sec_Network}, we provide a definition for overlap times in the single station setting and show how to extend overlap times to tandem networks of queues.  In Section \ref{same_cust}, we derive the joint distribution of overlap times of the same customer pair in a two dimensional tandem queue.  We use the joint distribution to calculate the mean, variance, and covariance of the overlaps in the tandem queue.  In Section \ref{diff_cust}, we analyze the setting where the customer pair in each station is different and also compute the joint tail distribution.  Finally, numerical examples are provided in Section \ref{numerics} and a conclusion is provided in Section \ref{conclusion}.  All the major proofs are provided in the Appendix or Section \ref{appendix}.

\section{Overlap Times in Networks of Queues} \label{sec_Network}

In this section, we study the infinite server queue with the intention of understanding how much time do adjacent customers spend in the system together.  A similar type of analysis has been completed \citet{kang2021queueing, palomooverlap}, however, in this work we extend the analysis to multiple queues.  Before we get into studying overlap times in the network setting, we first review overlap times in the single station case.  The review of single station case will be very helpful in the analysis that follows in the sequel.  

\subsection{Single Station Overlap Time Distribution}

In this section, we consider the GI/GI/$\infty$ queue.  Let $A_i$ be the arrival time of the $i^{th}$ customer, the inter-arrival time between the $i^{th}$ and $(i+1)^{th}$ customers is given by $A_{i+1} - A_i$.  In the renewal case, the inter-arrival times are i.i.d with cumulative distribution function (cdf) $F(x)$.  We also assume that $S_i$ be the service time of the $i^{th}$ customer and the service times are i.i.d with cdf $G(x)$.   Following \citet{palomooverlap}, the departure time for the $n^{th}$ customer is given by 
\begin{eqnarray}
D_n &=& S_n + A_{n}.
\end{eqnarray}
Given the arrival and departure time for each customer, it is possible to describe the overlap time between any pair of customers.  The overlap time between the $n^{th}$ and $(n+k)^{th}$ customers is given by
\begin{eqnarray}
O_{n,n+k} &=& \left( \min( D_n , D_{n+k} ) - A_{n+k} \right)^+ \\
&=& \left( \min( A_n + S_n , A_{n+k} + S_{n+k} ) - A_{n+k} \right)^+ \\
&=&  \left( (D_n - A_{n+k} )^+ \wedge S_{n+k} \right) \\
&=&  \left(  (S_n - \left( A_{n+k} - A_{n} \right) )^+ \wedge S_{n+k} \right) .
\end{eqnarray}

With this representation of the overlap time in terms of the service times and difference in arrival times, \citet{palomo2021overlap} were able to prove a result that describes the tail distribution of the overlap time.  Note that when the arrival and service distributions are independent, the overlap time distribution does not depend on $n$ and only depends on the gap between the two customers, which in this case is $k$.  In the sequel, we will use this representation to generalize it to the network setting and derive new results regarding the tail joint distribution of the overlap times of two stations.  

\subsection{Overlap in a Tandem Network}

Now that we understand the overlap times in the single station infinite server queue, we now describe the overlap times in a tandem network queue. Let $O^{(i)}_{n,n+k}$ be the overlap between the $n^{th}$ and $(n+k)^{th}$ customer in the $i^{th}$ station, then $O^{(i)}_{n,n+k}$ can be written as
\begin{eqnarray}
O^{(i)}_{n,n+k} &=& \left( \min( D^{(i)}_n , D^{(i)}_{n+k} ) - \max \left(D^{(i-1)}_n , D^{(i-1)}_{n+k} \right) \right)^+ \quad i \geq 1 
\end{eqnarray}
where $D^{0}_{n} = A_{n}$ and $D^{0}_{n+k} = A_{n+k}$.  Thus, the overlap between customer $n$ and $n+k$ in the $i^{th}$ station is equal to the difference of the \emph{first} customer to depart from station $i$ and the \emph{last} customer to depart station $i-1$.  

It is important to condition on whether the $n^{th}$ customer will leave before the $(n+k)^{th}$ customer.  Thus, we let $\mathcal{E}_{(\ell)}$ be the set corresponding to the event $A_{n+k} - A_n  \geq  \sum^{\ell -1}_{j=1} \left(S^{(j)}_n - S^{(j)}_{n+k} \right)$ or the event $D^{(\ell-1)}_{n+k} \geq D^{(\ell-1)}_{n} $.  This implies that the $n^{th}$ customer leaves the $(\ell -1)^{th}$ queue before the $(n+k)^{th}$ customer.  One should also note that 

\begin{eqnarray*}
     O^{(\ell)}_{n,n+k} 
    &=& \begin{cases}
    \min \left( S^{(1)}_n + \sum^{\ell -1}_{j=2} \left( S^{(j)}_n - S^{(j-1)}_{n+k} \right) + A_n - A_{n+k} , S^{(\ell)}_{n+k} \right)^+ , \mathrm{if} \ \mathcal{E}_{(\ell)} \ \mathrm{holds} \\
    \min \left( S^{(\ell)}_n , A_{n+k} - A_n + S^{(\ell)}_{n+k} + \sum^{\ell -1}_{j=2} \left(S^{(j)}_{n+k}  - S^{(j-1)}_n \right) \right)^+, \mathrm{if} \ \mathcal{E}^c_{(\ell)} \ \mathrm{holds}.
    \end{cases}
\end{eqnarray*}

Now we can use the expression for the overlap time of the two dimensional queueing network to compute both the joint and marginal distributions of the overlap time at each station.

\section{Two Stations with the Same Customer Pair} \label{same_cust}

In this section, we consider the overlap times of the $(M/M/\infty)^2$ queue where both stations have the same service time distribution.  We are particularly interested in computing the joint tail distribution of the overlap times of the same customer pair i.e. $$\mathbb{P} \left( O^{(1)}_{n,n+k} > x, O^{(2)}_{n,n+k} > y  \right). $$  In what follows, we will compute the tail joint distribution and the marginal distributions for each station.  We start with the joint distribution first.  

\subsection{The Joint Distribution}

In this section, we derive the joint tail distribution for the overlap times of the $(M/M/\infty)^2$ tandem queue.  

\begin{theorem}  \label{joint1}
Let $O^{(1)}_{n,n+k}$ and $O^{(2)}_{n,n+k}$ be the overlap times for the $n^{th}$ and $(n+k)^{th}$ customers in the $(M/M/\infty)^2$ tandem queue with identically distributed stations, then the joint distribution of the overlap times is equal to
\begin{eqnarray*}
\mathbb{P} \left( O^{(1)}_{n,n+k} > x, O^{(2)}_{n,n+k} > y  \right) &=& \frac{e^{-2\mu (x+y)} }{2 } \left( \frac{\lambda}{\lambda+\mu} \right)^{k}.
\end{eqnarray*}

\begin{proof}
The first step in the proof is to decompose the joint probability into two different components that separate the space into which queue departs out of the second to last queue first.  This yields the following expression.   
\begin{eqnarray*}
\lefteqn{ \mathbb{P} \left( O^{(1)}_{n,n+k} > x, O^{(2)}_{n,n+k} > y  \right) } \\  &=& \mathbb{P} \left( O^{(1)}_{n,n+k} > x, O^{(2)}_{n,n+k} > y, \mathcal{E}_{(2)}   \right) + \mathbb{P} \left( O^{(1)}_{n,n+k} > x, O^{(2)}_{n,n+k} > y, \mathcal{E}^c_{(2)}   \right) \\
&=& \mathbb{P} \left( S^{(1)}_n - \mathcal{A}_{k} > x , S^{(1)}_{n+k} > x, S^{(2)}_{n+k} > y, S^{(1)}_{n} + S^{(2)}_{n} - \mathcal{A}_{k} - S^{(1)}_{n+k} > y , \mathcal{E}_{(2)} \right) \\
&+& \mathbb{P} \left( S^{(1)}_n -  \mathcal{A}_{k} > x , S^{(1)}_{n+k} > x, S^{(2)}_{n} > y, S^{(1)}_{n+k} + S^{(2)}_{n+k} + \mathcal{A}_{k} - S^{(1)}_{n} > y,  \mathcal{E}^c_{(2)} \right).
\end{eqnarray*}
Now it remains to compute the last two probabilities, which can be done by explicit integration of two five dimensional integrals.  The first of two probability terms is given by 
\begin{eqnarray*}
\lefteqn{ \mathbb{P} \left( S^{(1)}_n - \mathcal{A}_{k} > x , S^{(1)}_{n+k} > x, S^{(2)}_{n+k} > y, S^{(1)}_{n} + S^{(2)}_{n} - \mathcal{A}_{k} - S^{(1)}_{n+k} > y  , \mathcal{E}_{(2)} \right)  } \\ 
&=& \int^{}_{\mathbb{R}^5_+}  \{ z_1 - z_5 > x \} \cdot \{ z_3 > x \} \cdot \{ z_4 > y  \} \cdot \{ z_1 + z_2 - z_5 - z_3 > y \} \cdot \{ z_5 + z_3 > z_1 \} \\
&\times& \frac{\lambda^k z_5^{k-1}}{\Gamma(k)} e^{-\lambda z_5} \mu^4 e^{-\mu (z_1 + z_2 + z_3 + z_4) } dz_1 dz_2 dz_3 dz_4 dz_5 \\
&=& e^{-\mu y} \int^{}_{\mathbb{R}^4_+}   \{ z_1 - z_5 > x \} \cdot \{ z_3 > x \}  \cdot \{ z_2  > y + z_5 + z_3 - z_1 \} \cdot \{ z_5 + z_3 > z_1 \} \\
&\times& \frac{\lambda^k z_5^{k-1}}{\Gamma(k)} e^{-\lambda z_5} \mu^3 e^{-\mu (z_1 + z_2 + z_3 ) } dz_1 dz_2 dz_3 dz_5 \\
&=& e^{-\mu y} \int^{}_{\mathbb{R}^3_+}   \{ z_1 > x + z_5 \} \cdot \{ z_3 > x \}  \cdot \{ z_5 + z_3 > z_1 \} \\
&\times& \frac{\lambda^k z_5^{k-1}}{\Gamma(k)} e^{-\lambda z_5} \mu^2 e^{-\mu (z_1 + z_3 ) } e^{-\mu (y + z_5 + z_3 - z_1)^+ } dz_1 dz_3 dz_5 \\
&=& e^{-\mu y} \int^{\infty}_{0} \int^{\infty}_{0} \int^{z_5 + z_3}_{x + z_5} \{ z_3 > x \}   \frac{\lambda^k z_5^{k-1}}{\Gamma(k)} e^{-\lambda z_5} \mu^2 e^{-\mu (z_1 + z_3 ) } e^{-\mu (y + z_5 + z_3 - z_1)^+ } dz_1 dz_3 dz_5 \\
&=& e^{-\mu y} \int^{\infty}_{0} \int^{\infty}_{0} \int^{z_5 + z_3}_{x + z_5} \{ z_3 > x \}   \frac{\lambda^k z_5^{k-1}}{\Gamma(k)} e^{-\lambda z_5} \mu^2 e^{-\mu (z_1 + z_3 ) } e^{-\mu (y + z_5 + z_3 - z_1) } dz_1 dz_3 dz_5 \\
&=& e^{-2\mu y} \int^{\infty}_{0} \int^{\infty}_{0} \int^{z_5 + z_3}_{x + z_5} \{ z_3 > x \}   \frac{\lambda^k z_5^{k-1}}{\Gamma(k)} e^{-\lambda z_5} \mu^2 e^{-\mu z_3 } e^{-\mu ( z_5 + z_3) } dz_1 dz_3 dz_5 \\
&=& e^{-2\mu y} \int^{\infty}_{0} \int^{\infty}_{0} \{ z_3 > x \}   \frac{\lambda^k z_5^{k-1}}{\Gamma(k)} e^{-\lambda z_5} \mu^2 e^{-\mu z_3 } e^{-\mu ( z_5 + z_3) } \left( z_3 - x \right) dz_3 dz_5 \\
&=&  \frac{e^{-2\mu (x+y)}}{4 } \int^{\infty}_{0}  \frac{\lambda^k z_5^{k-1}}{\Gamma(k)} e^{-(\lambda+\mu) z_5}  dz_5 \\
&=&  \frac{e^{-2\mu (x+y)}}{4 } \left( \frac{\lambda}{\lambda+\mu} \right)^k .
\end{eqnarray*}

Finally, the second probability term yields
\begin{eqnarray*}
\lefteqn{ \mathbb{P} \left( S^{(1)}_n -  \mathcal{A}_{k} > x , S^{(1)}_{n+k} > x, S^{(2)}_{n} > y, S^{(1)}_{n+k} + S^{(2)}_{n+k} + \mathcal{A}_{k} - S^{(1)}_{n} > y  , \mathcal{E}^c_{(2)} \right)  } \\ 
&=&  \int^{}_{\mathbb{R}_+^5}  \{ z_1 > z_5 + (z_3 \vee x) \} \cdot \{ z_3 > x \} \cdot \{ z_2 > y  \} \cdot \{ z_3 + z_4 + z_5 - z_1 > y \} \\
&\times& \frac{\lambda^k z_5^{k-1}}{\Gamma(k)} e^{-\lambda z_5} \mu^4 e^{-\mu (z_1 + z_2 + z_3 + z_4) } dz_1 dz_2 dz_3 dz_4 dz_5 \\
&=& e^{-\mu y}  \int^{}_{\mathbb{R}_+^4}    \{ z_1 > z_5 + (z_3 \vee x) \} \cdot \{ z_3 > x \} \cdot  \{ z_4 > y - z_3 - z_5 + z_1 \} \\
&\times& \frac{\lambda^k z_5^{k-1}}{\Gamma(k)} e^{-\lambda z_5} \mu^3 e^{-\mu (z_1 + z_3 + z_4) } dz_1 dz_3 dz_4 dz_5 \\
&=& e^{-\mu y}  \int^{}_{\mathbb{R}_+^3}    \{ z_1 > z_5 + (z_3 \vee x) \}\cdot \{ z_3 > x \} \\
&\times& \frac{\lambda^k z_5^{k-1}}{\Gamma(k)} e^{-\lambda z_5} \mu^2 e^{-\mu (z_1 + z_3 + z_4) } e^{-\mu (y - z_3 - z_5 + z_1)^+ } dz_1 dz_3 dz_5 \\ 
&=& e^{-2\mu y} \int^{\infty}_{0} \int^{\infty}_{0} \int^{\infty}_{z_5 + (z_3 \vee x) }   \{ z_3 > x \} \frac{\lambda^k z_5^{k-1}}{\Gamma(k)} e^{- (\lambda-\mu) z_5} \mu^2 e^{-2\mu z_1  } dz_1 dz_3 dz_5 \\
&=& e^{-2\mu y} \int^{\infty}_{0} \int^{\infty}_{x} \frac{\lambda^k z_5^{k-1}}{\Gamma(k)} e^{- (\lambda+\mu) z_5} \frac{\mu}{2} e^{-2\mu z_3}  dz_3 dz_5 \\
&=& e^{-2\mu (x+y)} \int^{\infty}_{0}  \frac{\lambda^k z_5^{k-1}}{\Gamma(k)} e^{- (\lambda+\mu) z_5} \frac{1}{4}  dz_5 \\
&=& \frac{e^{-2\mu (x+y)} \lambda^k }{4 \Gamma(k)} \int^{\infty}_{0}  z_5^{k-1} e^{- (\lambda+\mu) z_5} dz_5 \\
&=& \frac{e^{-2\mu (x+y)} }{4 } \left( \frac{\lambda}{\lambda+\mu} \right)^k  .
\end{eqnarray*}

Now by adding the two probability terms together, we obtain the result and complete the proof. 
\end{proof}
\end{theorem}

\subsection{The Marginal Distribution}

In this section, we compute the marginal tail distribution of the overlap time in each station.  We assume that all the stations have the same service distributions and each is exponential with rate $\mu$.   In order to compute the tail probability, we need a series of lemmas, which are essential for deriving the main results.

\begin{theorem} \label{marginals}
Let $O^{(1)}_{n,n+k}$ and $O^{(2)}_{n,n+k}$ be the overlap times for the $n^{th}$ and $(n+k)^{th}$ customers in the $(M/M/\infty)^2$ tandem queue with identically distributed stations, then the marginal distributions of the overlap times is equal to
\begin{eqnarray*}
\mathbb{P} \left(  O^{(1)}_{n,n+k} > x  \right) &=&  e^{-2\mu x} \alpha^k  \\
\mathbb{P} \left(  O^{(2)}_{n,n+k} > y  \right) &=& \frac{e^{-2\mu y} }{2 } \alpha^k  + \frac{e^{-2\mu y} }{2 } \cdot   \alpha^k   \cdot \left( \frac{\mu}{\lambda + \mu} \right) \cdot k . 
\end{eqnarray*}
Moreover, the remaining joint probabilities have the following expressions
\begin{eqnarray*}
\mathbb{P} \left(  O^{(1)}_{n,n+k} \leq x , O^{(2)}_{n,n+k} \leq y  \right) &=& 1  - \left( \frac{e^{-2\mu y} }{2 } + \frac{e^{-2\mu y} }{2 } \cdot \left( \frac{\mu k }{\lambda + \mu} \right)  - \frac{e^{-2\mu (x+y)}}{2} +  e^{-2\mu x} \right)\alpha^k   \\
\mathbb{P} \left(  O^{(1)}_{n,n+k} \leq x , O^{(2)}_{n,n+k} > y  \right) &=& \left( \frac{e^{-2\mu y} }{2 } + \frac{e^{-2\mu y} }{2 } \cdot \left( \frac{\mu}{\lambda + \mu} \right) \cdot k - e^{-2\mu (x+y)} \right) \alpha^k  \\
\mathbb{P} \left(  O^{(1)}_{n,n+k} > x , O^{(2)}_{n,n+k} \leq y  \right) &=&  \left( e^{-2\mu x} -  e^{-2\mu (x+y)} \right) \alpha^k.
\end{eqnarray*}

\begin{proof}
\begin{eqnarray*}
\mathbb{P} \left(  O^{(2)}_{n,n+k} > y  \right) &=&  \mathbb{P} \left( O^{(1)}_{n,n+k} > 0, O^{(2)}_{n,n+k} > y  \right)  + \mathbb{P} \left( O^{(1)}_{n,n+k} = 0, O^{(2)}_{n,n+k} > y  \right) \\
&=&\frac{e^{-2\mu y} }{2 } \alpha^k  + \mathbb{P} \left( O^{(1)}_{n,n+k} = 0, O^{(2)}_{n,n+k} > y  \right).
\end{eqnarray*}
Now it remains to compute $\mathbb{P} \left( O^{(1)}_{n,n+k} = 0, O^{(2)}_{n,n+k} > y  \right)$.

\begin{eqnarray*}
\lefteqn{ \mathbb{P} \left( O^{(1)}_{n,n+k} = 0, O^{(2)}_{n,n+k} > y  \right) } \\ &=&  \int^{}_{\mathbb{R}^5_+}  \{ z_1 - z_5 < 0 \} \cdot \{ z_3 > 0 \} \cdot \{ z_4 > y  \} \cdot \{ z_1 + z_2 - z_5 - z_3 > y \}  \\
&\times& \frac{\lambda^k z_5^{k-1}}{\Gamma(k)} e^{-\lambda z_5} \mu^4 e^{-\mu (z_1 + z_2 + z_3 + z_4) } dz_1 dz_2 dz_3 dz_4 dz_5 \\
&=& e^{-\mu y} \int^{}_{\mathbb{R}^4_+}  \{ z_1  < z_5 \} \cdot  \{ z_2 > y +z_3 + z_5 - z_1 \}  \\
&\times& \frac{\lambda^k z_5^{k-1}}{\Gamma(k)} e^{-\lambda z_5} \mu^3 e^{-\mu (z_1 + z_2 + z_3 ) } dz_1 dz_2 dz_3  dz_5 \\
&=& e^{-2\mu y} \int^{}_{\mathbb{R}^3_+}  \{ z_1  < z_5 \}  \frac{\lambda^k z_5^{k-1}}{\Gamma(k)} e^{-(\lambda+\mu) z_5} \mu^2 e^{-2 \mu z_3  } dz_1 dz_3  dz_5 \\
&=& e^{-2\mu y} \int^{}_{\mathbb{R}^2_+} \frac{\lambda^k z_5^k}{\Gamma(k)} e^{-(\lambda+\mu) z_5} \mu^2 e^{-2 \mu z_3  } dz_3  dz_5 \\
&=& e^{-2\mu y} \int^{}_{\mathbb{R}^1_+} \frac{\lambda^k z_5^k}{\Gamma(k)} e^{-(\lambda+\mu) z_5} \frac{\mu}{2} dz_5 \\
&=& \frac{\mu e^{-2\mu y}}{2} \frac{\lambda^k}{(\lambda +\mu)^{k+1}} \frac{\Gamma(k+1)}{\Gamma(k)} \\
&=& \frac{\mu k e^{-2\mu y}}{2} \frac{\lambda^k}{(\lambda +\mu)^{k+1}} .
\end{eqnarray*}
The remaining terms can be determined by the marginals and the tail distribution using the following expressions
\begin{eqnarray*}
\mathbb{P} \left( O^{(1)}_{n,n+k} \leq x, O^{(2)}_{n,n+k} \leq y \right) 
&=& \mathbb{P} \left( O^{(1)}_{n,n+k} \leq x \right) - \mathbb{P} \left( O^{(2)}_{n,n+k} > y \right) \\
&+& \mathbb{P} \left( O^{(1)}_{n,n+k} > x, O^{(2)}_{n,n+k} > y \right) 
\end{eqnarray*}
\begin{eqnarray*}
\mathbb{P} \left( O^{(1)}_{n,n+k} \leq x, O^{(2)}_{n,n+k} > y \right) &=&  \mathbb{P} \left( O^{(2)}_{n,n+k} > y  \right) - \mathbb{P} \left( O^{(1)}_{n,n+k} > x, O^{(2)}_{n,n+k} > y \right) 
\end{eqnarray*}
\begin{eqnarray*}
\mathbb{P} \left( O^{(1)}_{n,n+k} > x, O^{(2)}_{n,n+k} \leq y \right) &=&  \mathbb{P} \left( O^{(1)}_{n,n+k} > x  \right) - \mathbb{P} \left( O^{(1)}_{n,n+k} > x, O^{(2)}_{n,n+k} > y \right) .
\end{eqnarray*}
By substituting the expressions in for each of the probabilities completes the proof.
\end{proof}
\end{theorem}


\subsection{Tail Distribution of Sum of Overlap Times}

In this section, we compute the tail distribution of the sum of the overlap times in both stations.  The sum is of interest as it is important to know how long any two customers will overlap in the full system.  This is challenging as they are dependent random variables, however, we provide a way of deriving the tail distribution of the sum in what follows.  

\begin{theorem}
Let $O^{(1)}_{n,n+k}$ and $O^{(2)}_{n,n+k}$ be the overlap times for the $n^{th}$ and $(n+k)^{th}$ customers in the $(M/M/\infty)^2$ tandem queue with identically distributed stations, then the distribution of the sum of the overlap times is equal to
\begin{eqnarray*}
\mathbb{P} \left(  O^{(1)}_{n,n+k} + O^{(2)}_{n,n+k} > \ell  \right)  &=&  \alpha^k  (1+2 \ell \mu) e^{-2 \ell \mu} +  e^{-2 \mu \ell}  \frac{\alpha^k}{2} (1-\alpha) k  +  e^{-2 \mu \ell}  \frac{\alpha^k}{2}.
\end{eqnarray*}
\begin{proof}
We start the proof by letting $\alpha = \frac{\lambda}{\lambda+\mu}$.  Next, we take the partial derivative of the tail distribution with respect to both $x$ and $y$.  This yields the following joint density of the overlap times in both stations 
\begin{align*}
    \frac{\partial^2 \mathbb{P} \left( O^{(1)}_{n,n+k} > x, O^{(2)}_{n,n+k} > y  \right) }{\partial y \partial x} &= \frac{\partial}{\partial x} \left( \frac{\partial \mathbb{P} \left( O^{(1)}_{n,n+k} > x, O^{(2)}_{n,n+k} > y  \right) }{\partial y} \right)\\
    &= 4 \mu^2 e^{-2\mu (x+y)} \alpha^k.
\end{align*}

So,
\begin{align*}
    \mathbb{P} \left(  O^{(1)}_{n,n+k} + O^{(2)}_{n,n+k} > \ell  \right) &= \int_0^\infty \int_{(\ell-y)^+}^\infty 4 \mu^2 e^{-2\mu (x+y)} \alpha^k dx dy.
\end{align*}
But, since we also have the constraint $x \in [0,\infty]$, we modify the expression to be:
\begin{eqnarray*}
    \mathbb{P} \left(  O^{(1)}_{n,n+k} + O^{(2)}_{n,n+k} > \ell  \right) &=& \int_0^l \int_{l-y}^\infty 4 \mu^2 e^{-2\mu (x+y)} \alpha^k dx dy \\
    &+& \int_l^\infty \int_0^\infty 4 \mu^2 e^{-2\mu (x+y)} \alpha^k dx dy \\
    &=& 2 \ell \mu \alpha^k e^{-2 \ell \mu} + \alpha^k e^{-2 \ell \mu} \\
    &=& \alpha^k  (1+2 \ell \mu) e^{-2 \ell \mu} .
\end{eqnarray*}

\begin{eqnarray*}
    \lefteqn{ \mathbb{P} \left(  O^{(1)}_{n,n+k} + O^{(2)}_{n,n+k} > \ell  \right) } \\ 
    &=&  \mathbb{P} \left(  O^{(1)}_{n,n+k} + O^{(2)}_{n,n+k} > \ell | [ O^{(1)}_{n,n+k} > 0 , O^{(2)}_{n,n+k} > 0 ]  \right)  \\
    &+&  \mathbb{P} \left(  O^{(1)}_{n,n+k} + O^{(2)}_{n,n+k} > \ell | [ O^{(1)}_{n,n+k} = 0  , O^{(2)}_{n,n+k} > 0 ]   \right)  \\
    &+&  \mathbb{P} \left(  O^{(1)}_{n,n+k} + O^{(2)}_{n,n+k} > \ell | [O^{(1)}_{n,n+k} > 0  , O^{(2)}_{n,n+k} = 0]   \right)  \\
    &+&  \mathbb{P} \left(  O^{(1)}_{n,n+k} + O^{(2)}_{n,n+k} > \ell | [O^{(1)}_{n,n+k} = 0  , O^{(2)}_{n,n+k} = 0]   \right)  \\
    &=&  \mathbb{P} \left(  O^{(1)}_{n,n+k} + O^{(2)}_{n,n+k} > \ell | [O^{(1)}_{n,n+k} > 0 , O^{(2)}_{n,n+k} > 0]   \right)  \\
    &+&  \mathbb{P} \left(   O^{(2)}_{n,n+k} > \ell , O^{(1)}_{n,n+k} = 0   \right) + \mathbb{P} \left(   O^{(1)}_{n,n+k} > \ell , O^{(2)}_{n,n+k} = 0   \right) \\
    &=& \alpha^k  (1+2 \ell \mu) e^{-2 \mu \ell}  +   \left( \mathbb{P} \left(  O^{(2)}_{n,n+k} > \ell   \right) - \mathbb{P} \left(  O^{(1)}_{n,n+k} > 0 , O^{(2)}_{n,n+k} > \ell   \right) \right) \\
    &+&   \left( \mathbb{P} \left(  O^{(1)}_{n,n+k} > \ell   \right) - \mathbb{P} \left(  O^{(1)}_{n,n+k} > \ell , O^{(2)}_{n,n+k} > 0   \right)\right) \\
    &=& \alpha^k  (1+2 \ell \mu) e^{-2 \ell \mu} + \left( e^{-2 \mu \ell}  \frac{\alpha^k}{2} + e^{-2 \mu \ell}  \frac{\alpha^k}{2} (1-\alpha) k - e^{-2 \mu \ell}  \frac{\alpha^k}{2} \right) \\
    &+&   \left(   e^{-2 \mu \ell}  \alpha^k  - e^{-2 \mu \ell}  \frac{\alpha^k}{2} \right) \\
    &=& \alpha^k  (1+2 \ell \mu) e^{-2 \ell \mu} +  e^{-2 \mu \ell}  \frac{\alpha^k}{2} (1-\alpha) k  +  e^{-2 \mu \ell}  \frac{\alpha^k}{2}.
\end{eqnarray*}
This completes the proof.
\end{proof}
\end{theorem}

\subsection{Moments and Covariance}

\begin{proposition}
Let $O^{(1)}_{n,n+k}$ and $O^{(2)}_{n,n+k}$ be the overlap times for the $n^{th}$ and $(n+k)^{th}$ customers in the $(M/M/\infty)^2$ tandem queue with identically distributed stations, then the mean, variance and covariance of the two overlap times are equal to
\begin{eqnarray}
\mathbb{E}\left[ O^{(1)}_{n,n+k}  \right] &=&  \left( \frac{\lambda}{\lambda+\mu} \right)^k  \frac{1}{2 \mu }   \\
\mathbb{E}\left[ O^{(2)}_{n,n+k}   \right] &=&  \frac{1 }{4\mu } \left( \frac{\lambda}{\lambda+\mu} \right)^{k} +  \frac{k }{4 } \left( \frac{\lambda}{\lambda+\mu} \right)^{k} \cdot \left( \frac{1}{\lambda + \mu} \right)    \\
\mathrm{Var}\left[ O^{(1)}_{n,n+k}  \right] &=&  \left( \frac{\lambda}{\lambda+\mu} \right)^k \frac{1}{ 2\mu^2 }  - \left( \frac{\lambda}{\lambda+\mu} \right)^{2k} \frac{1}{4 \mu^2 } \\
\mathrm{Var}\left[ O^{(2)}_{n,n+k}  \right] &=&  \left( \frac{\lambda}{\lambda+\mu} \right)^k \frac{1}{ 4\mu^2 }  + \frac{1}{ 4\mu^2 } \left( \left( \frac{\lambda}{\lambda+\mu} \right)^{k}  \cdot \left( \frac{\mu}{\lambda + \mu} \right) \cdot k  \right) \nonumber \\
&-& \left( \frac{1 }{4\mu } \left( \frac{\lambda}{\lambda+\mu} \right)^{k} +  \frac{k }{4 } \left( \frac{\lambda}{\lambda+\mu} \right)^{k} \cdot \left( \frac{1}{\lambda + \mu} \right) \right)^2 \\
\mathrm{Cov}\left[ O^{(1)}_{n,n+k} , O^{(2)}_{n,n+k}   \right] &=&  \frac{\left( \frac{\lambda}{\lambda+\mu} \right)^k}{8 \mu^2} - \frac{\left( \frac{\lambda}{\lambda+\mu} \right)^{2k}}{8 \mu^2} \left(  1 + \frac{k \mu}{\lambda + \mu} \right) .
\end{eqnarray}

\end{proposition}
\begin{proof}
The mean of $O^{(1)}_{n,n+k}$ can be obtained by integrating the tail distribution.  This is equivalent to
\begin{eqnarray*}
\mathbb{E}\left[ O^{(1)}_{n,n+k}   \right] &=& \int^{\infty}_{0}  \mathbb{P} \left( O^{(1)}_{n,n+k} > x \right) dx   \\
&=& \left( \frac{\lambda}{\lambda+\mu} \right)^k \int^{\infty}_{0} e^{-2\mu x} dx   \\
&=& \left( \frac{\lambda}{\lambda+\mu} \right)^k \frac{1}{2 \mu }    .
\end{eqnarray*}
The mean of $O^{(2)}_{n,n+k}$ can also be obtained by integrating the tail distribution.  This is equivalent to
\begin{eqnarray*}
\mathbb{E}\left[ O^{(2)}_{n,n+k}   \right] &=& \int^{\infty}_{0}  \mathbb{P} \left( O^{(2)}_{n,n+k} > x \right) dx   \\
&=&  \int^{\infty}_{0} \left( \frac{e^{-2\mu x} }{2 } \left( \frac{\lambda}{\lambda+\mu} \right)^{k} + \frac{e^{-2\mu x} }{2 } \cdot   \left( \frac{\lambda}{\lambda+\mu} \right)^{k}  \cdot \left( \frac{\mu}{\lambda + \mu} \right) \cdot k \right) dx   \\
&=&  \frac{1 }{4\mu } \left( \frac{\lambda}{\lambda+\mu} \right)^{k} +  \frac{k }{4 } \left( \frac{\lambda}{\lambda+\mu} \right)^{k} \cdot \left( \frac{1}{\lambda + \mu} \right) .
\end{eqnarray*}
The second moment of the first overlap time is equal to 
\begin{eqnarray*}
\mathbb{E} \left[ (O^{(1)}_{n,n+k} )^2   \right] &=&\int^{\infty}_{0} 2 x  \mathbb{P} \left( O^{(1)}_{n,n+k} > x \right) dx   \\
&=& \left( \frac{\lambda}{\lambda+\mu} \right)^k \int^{\infty}_{0} 2 x e^{-2\mu x} dx   \\
&=& \left( \frac{\lambda}{\lambda+\mu} \right)^k \frac{1}{ 2\mu^2 }   .
\end{eqnarray*}
The second moment of the second overlap time is equal to 
\begin{eqnarray*}
\lefteqn{ \mathbb{E}\left[ (O^{(2)}_{n,n+k})^2   \right] } \\ &=& \int^{\infty}_{0} 2 x  \mathbb{P} \left( O^{(2)}_{n,n+k} > x \right) dx   \\
&=& \int^{\infty}_{0} 2 x \left( \frac{e^{-2\mu x} }{2 } \left( \frac{\lambda}{\lambda+\mu} \right)^{k} + \frac{e^{-2\mu x} }{2 } \cdot   \left( \frac{\lambda}{\lambda+\mu} \right)^{k}  \cdot \left( \frac{\mu}{\lambda + \mu} \right) \cdot k  \right) dx  \\
&=& \left( \frac{\lambda}{\lambda+\mu} \right)^k \frac{1}{ 4\mu^2 }  + \frac{1}{ 4\mu^2 } \left( \left( \frac{\lambda}{\lambda+\mu} \right)^{k}  \cdot \left( \frac{\mu}{\lambda + \mu} \right) \cdot k  \right)  .
\end{eqnarray*}
\begin{eqnarray*}
\mathrm{Var}\left[ O^{(1)}_{n,n+k}  \right] &=&  \left( \frac{\lambda}{\lambda+\mu} \right)^k \frac{1}{ 2\mu^2 }  - \left( \frac{\lambda}{\lambda+\mu} \right)^{2k} \frac{1}{4 \mu^2 } 
\end{eqnarray*}
\begin{eqnarray}
\mathrm{Var}\left[ O^{(2)}_{n,n+k}  \right] &=&  \left( \frac{\lambda}{\lambda+\mu} \right)^k \frac{1}{ 4\mu^2 }  + \frac{1}{ 4\mu^2 } \left( \left( \frac{\lambda}{\lambda+\mu} \right)^{k}  \cdot \left( \frac{\mu}{\lambda + \mu} \right) \cdot k  \right) \\
&-& \left( \frac{1 }{4\mu } \left( \frac{\lambda}{\lambda+\mu} \right)^{k} +  \frac{k }{4 } \left( \frac{\lambda}{\lambda+\mu} \right)^{k} \cdot \left( \frac{1}{\lambda + \mu} \right) \right)^2 .
\end{eqnarray}
Finally, we calculate the product of the two overlap times as
\begin{eqnarray*}
\mathbb{E}\left[ O^{(1)}_{n,n+k} \cdot O^{(2)}_{n,n+k}   \right] &=& \int^{\infty}_{0} \int^{\infty}_{0} \mathbb{P} \left( O^{(1)}_{n,n+k} > x, O^{(2)}_{n,n+k} > y  \right) dx dy  \\
&=& \int^{\infty}_{0} \int^{\infty}_{0} \frac{e^{-2\mu (x+y)} }{2 } \left( \frac{\lambda}{\lambda+\mu} \right)^k  dx dy  \\
&=& \left( \frac{\lambda}{\lambda+\mu} \right)^k \int^{\infty}_{0} \int^{\infty}_{0} \frac{e^{-2\mu (x+y)} }{2 }  dx dy  \\
&=& \left( \frac{\lambda}{\lambda+\mu} \right)^k \int^{\infty}_{0} \frac{e^{-2\mu y} }{4\mu } dy  \\
&=& \left( \frac{\lambda}{\lambda+\mu} \right)^k  \frac{1}{8 \mu^2 }   .
\end{eqnarray*}
\begin{eqnarray*}
\mathrm{Cov}\left[ O^{(1)}_{n,n+k} , O^{(2)}_{n,n+k}   \right] &=&  \frac{\left( \frac{\lambda}{\lambda+\mu} \right)^k}{8 \mu^2} - \frac{\left( \frac{\lambda}{\lambda+\mu} \right)^{2k}}{8 \mu^2} \left(  1 + \frac{k \mu}{\lambda + \mu} \right) . 
\end{eqnarray*}
\end{proof}

\section{Two Stations with Different Customer Pairs} \label{diff_cust}
So far, the paper has considered identical stations with the same pair of customers i.e. (n,n+k) in both stations. Now we explore the case where the customer pair of interest shifts. We now use nine cases to exhaust all possible customer pairs and compute the joint tail distribution respectively in each case.

\subsection{Joint Tail Distribution}
\label{ovlp_diff_customers}
\begin{theorem} \label{ovlp_diff_customers_prop}
Let $O^{(1)}_{n,n+j}$ and $O^{(2)}_{m,m+k}$ be the overlap time for the $n^{th}$ and $(n+j)^{th}$ customers in the $(M/M/\infty)$ tandem queue and with the overlap time for the $m^{th}$ and $(m+k)^{th}$ customers in another $(M/M/\infty)$ tandem queue, then the joint probabilities for different $n,m,j,k$ are:
\begin{enumerate}
    \item when $n = m$, $j = k$:
    \begin{eqnarray*}
    \mathbb{P} \left( O^{(1)}_{n,n+j} > x, O^{(2)}_{m,m+k} > y  \right) 
    &=&  \frac{e^{-2\mu (x+y)}}{2} \left( \frac{\lambda}{\lambda+\mu} \right)^k
    \end{eqnarray*}
    \item when $n = m$, $j < k$, let $I$ denote ${\mathds{1}_{\lambda<\mu}}$:
    \begin{enumerate}
        \item $\lambda > \mu$
            \begin{eqnarray*}
                \lefteqn{ \mathbb{P} \left( O^{(1)}_{n,n+j} > x, O^{(2)}_{n,n+k} > y  \right) } \\
                &=& \frac{1}{4} e^{-\mu (3x+2y)}\cdot \frac{\lambda^{j}}{(\mu+\lambda)^{j}} \cdot \frac{\lambda^{k-j}}{(\lambda - \mu)^{k-j}} \cdot \frac{1}{\Gamma(k-j)} \cdot \gamma (k-j,(\lambda-\mu)x)\\
&+& \frac{1}{4} e^{-\mu (x+2y)}\cdot \frac{\lambda^{j}}{(\mu+\lambda)^{j}} \cdot \frac{\lambda^{k-j}}{(\lambda +\mu)^{k-j}} \cdot \left(1- \frac{1}{\Gamma(k-j)} \cdot \gamma (k-j,(\lambda+\mu)x)\right) \\
&+& \frac{\mu}{2} e^{-\mu (x+2y)} \cdot \frac{\lambda^{j}}{(\mu+\lambda)^{j}} \cdot \frac{\lambda^{k-j}}{(\mu+\lambda)^{k-j}}\left(\frac{k-j}{\mu+\lambda}-x\right)\\
&-&\frac{\mu}{2} e^{-\mu (x+2y)} \cdot \frac{\lambda^{j}}{(\mu+\lambda)^{j}} \cdot \frac{\lambda^{k-j}}{(\mu+\lambda)^{k-j}} \left(\frac{\gamma\left(k-j+1,x(\lambda+\mu)\right)}{(\lambda+\mu)\Gamma(k-j)} -x \frac{\gamma\left(k-j,x(\lambda+\mu)\right)}{\Gamma(k-j)} \right)\\
&+& \frac{1}{4} e^{-\mu (3x+2y)} \left( \frac{\lambda}{\mu+\lambda} \right)^j \cdot \left( \frac{\lambda}{\lambda-\mu} \right)^{k-j} \frac{\gamma \big( k-j, (\lambda-\mu)x \big)}{\Gamma(k-j)} \cdot \left( 2\mu x +1 - \frac{2 \mu j}{\mu+\lambda} \right)\\
    &+& \frac{1}{4} e^{-\mu (x+2y)} \cdot \left( \frac{\lambda}{\mu+\lambda} \right)^k \cdot \left( 1- \frac{\gamma \big( k-j, (\lambda+\mu)x \big)}{\Gamma(k-j)} \right) \cdot \left( 1+ \frac{2\mu j}{\mu + \lambda} \right)\\
    &+&\frac{1}{4} e^{-\mu (x+2y)} \cdot \left( \frac{\lambda}{\mu+\lambda} \right)^k \cdot \frac{k-j}{\lambda+\mu} \cdot \left( 1- \frac{\gamma \big( k-j+1, (\lambda+\mu)x \big)}{\Gamma(k-j)} \right) \\
    &-& \frac{1}{2} \mu e^{-\mu(3x+2y)} \left( \frac{\lambda}{\mu+\lambda} \right)^j \cdot \frac{j}{\lambda+\mu} \cdot \left( \frac{\lambda}{\lambda-\mu} \right)^{k-j} \cdot \frac{\gamma \big( k-j, (\lambda-\mu)x \big)}{\Gamma(k-j)}\\
    &-&\frac{1}{2} \mu e^{-\mu(x+2y)} \left( \frac{\lambda}{\mu+\lambda} \right)^k \cdot \frac{j}{\lambda+\mu} \cdot \left(1- \frac{\gamma \big( k-j, (\lambda+\mu)x \big)}{\Gamma(k-j)} \right)\\
    &-& \frac{1}{2} \mu e^{-\mu(3x+2y)} \left( \frac{\lambda}{\mu+\lambda} \right)^j \cdot \frac{k-j}{\lambda-\mu} \cdot \left( \frac{\lambda}{\lambda-\mu} \right)^{k-j} \cdot \frac{\gamma \big( k-j+1, (\lambda-\mu)x \big)}{\Gamma(k-j+1)}\\
    &-&\frac{1}{2} \mu e^{-\mu(x+2y)} \left( \frac{\lambda}{\mu+\lambda} \right)^k \cdot \frac{k-j}{\lambda+\mu} \cdot \left(1- \frac{\gamma \big( k-j+1, (\lambda+\mu)x \big)}{\Gamma(k-j+1)} \right).\\
            \end{eqnarray*}
        \item $\lambda = \mu$
            \begin{eqnarray*}
                \lefteqn{ \mathbb{P} \left( O^{(1)}_{n,n+j} > x, O^{(2)}_{n,n+k} > y  \right) } \\
                &=& \frac{1}{4} e^{-\mu (3x+2y)} \left( \frac{\lambda}{\mu+\lambda} \right)^j \cdot \left( \frac{\lambda}{\lambda-\mu} \right)^{k-j} \frac{-\gamma \big( k-j, (-\lambda+\mu)x \big)}{\Gamma(k-j)} \cdot \left( 2\mu x +1 - \frac{2 \mu j}{\mu+\lambda} \right)\\
                &+& \frac{1}{4} e^{-\mu (x+2y)} \cdot \left( \frac{\lambda}{\mu+\lambda} \right)^k \cdot \left( 1- \frac{\gamma \big( k-j, (\lambda+\mu)x \big)}{\Gamma(k-j)} \right) \cdot \left( 1+ \frac{2\mu j}{\mu + \lambda} \right)\\
                &+&\frac{1}{4} e^{-\mu (x+2y)} \cdot \left( \frac{\lambda}{\mu+\lambda} \right)^k \cdot \frac{k-j}{\lambda+\mu} \cdot \left( 1- \frac{\gamma \big( k-j+1, (\lambda+\mu)x \big)}{\Gamma(k-j)} \right) \\
                &-& \frac{1}{2} \mu e^{-\mu(3x+2y)} \left( \frac{\lambda}{\mu+\lambda} \right)^j \cdot \frac{j}{\lambda+\mu} \cdot \left( \frac{\lambda}{\lambda-\mu} \right)^{k-j} \cdot \frac{-\gamma \big( k-j, (-\lambda+\mu)x \big)}{\Gamma(k-j)}\\
                &-&\frac{1}{2} \mu e^{-\mu(x+2y)} \left( \frac{\lambda}{\mu+\lambda} \right)^k \cdot \frac{j}{\lambda+\mu} \cdot \left(1- \frac{\gamma \big( k-j, (\lambda+\mu)x \big)}{\Gamma(k-j)} \right)\\
                &+& \frac{1}{2} \mu e^{-\mu(3x+2y)} \left( \frac{\lambda}{\mu+\lambda} \right)^j \cdot \frac{k-j}{\lambda-\mu} \cdot \left( \frac{\lambda}{\lambda-\mu} \right)^{k-j} \cdot \frac{\gamma \big( k-j+1, (-\lambda+\mu)x \big)}{\Gamma(k-j+1)}\\
                &-&\frac{1}{2} \mu e^{-\mu(x+2y)} \left( \frac{\lambda}{\mu+\lambda} \right)^k \cdot \frac{k-j}{\lambda+\mu} \cdot \left(1- \frac{\gamma \big( k-j+1, (\lambda-\mu)x \big)}{\Gamma(k-j+1)} \right)\\
            \end{eqnarray*}
        \item $\lambda < \mu$
            \begin{eqnarray*}
                \lefteqn{ \mathbb{P} \left( O^{(1)}_{n,n+j} > x, O^{(2)}_{n,n+k} > y  \right) } \\
                &=& \frac{1}{4} e^{-\mu (3x+2y)}\cdot \frac{\lambda^{j}}{(\mu+\lambda)^{j}} \frac{\lambda^{(k-j)}}{\Gamma(k-j) } \\
&\times& \left(  \sum_{i=0}^{k-j-1} (-1)^{i} \cdot \frac{e^{(\mu-\lambda)x}(k-j-1)!}{(k-j-1-i)!\cdot \mu-\lambda)^{i+1}} \cdot x^{k-j-1-i} - (-1)^{k-j-1}\frac{(k-j-1)!}{\mu-\lambda)^{k-j-1}} \right)\\
&+& \frac{1}{4} e^{-\mu (x+2y)}\cdot \frac{\lambda^{j}}{(\mu+\lambda)^{j}} \cdot \frac{\lambda^{k-j}}{(\lambda +\mu)^{k-j}} \cdot \left(1- \frac{1}{\Gamma(k-j)} \cdot \gamma (k-j,(\lambda+\mu)x)\right) \\
&+& \frac{\mu}{2} e^{-\mu (x+2y)} \cdot \frac{\lambda^{j}}{(\mu+\lambda)^{j}} \cdot \frac{\lambda^{k-j}}{(\mu+\lambda)^{k-j}}\left(\frac{k-j}{\mu+\lambda}-x\right)\\
&-&\frac{\mu}{2} e^{-\mu (x+2y)} \cdot \frac{\lambda^{j}}{(\mu+\lambda)^{j}} \cdot \frac{\lambda^{k-j}}{(\mu+\lambda)^{k-j}} \left(\frac{\gamma\left(k-j+1,x(\lambda+\mu)\right)}{(\lambda+\mu)\Gamma(k-j)} -x \frac{\gamma\left(k-j,x(\lambda+\mu)\right)}{\Gamma(k-j)} \right)\\
&+& \frac{1}{4} e^{-\mu (3x+2y)} \left( \frac{\lambda}{\mu+\lambda} \right)^j \cdot \left( \frac{\lambda}{\lambda-\mu} \right)^{k-j} \frac{-\gamma \big( k-j, (-\lambda+\mu)x \big)}{\Gamma(k-j)} \cdot \left( 2\mu x +1 - \frac{2 \mu j}{\mu+\lambda} \right)\\
    &+& \frac{1}{4} e^{-\mu (x+2y)} \cdot \left( \frac{\lambda}{\mu+\lambda} \right)^k \cdot \left( 1- \frac{\gamma \big( k-j, (\lambda+\mu)x \big)}{\Gamma(k-j)} \right) \cdot \left( 1+ \frac{2\mu j}{\mu + \lambda} \right)\\
    &+&\frac{1}{4} e^{-\mu (x+2y)} \cdot \left( \frac{\lambda}{\mu+\lambda} \right)^k \cdot \frac{k-j}{\lambda+\mu} \cdot \left( 1- \frac{\gamma \big( k-j+1, (\lambda+\mu)x \big)}{\Gamma(k-j)} \right) \\
    &-& \frac{1}{2} \mu e^{-\mu(3x+2y)} \left( \frac{\lambda}{\mu+\lambda} \right)^j \cdot \frac{j}{\lambda+\mu} \cdot \left( \frac{\lambda}{\lambda-\mu} \right)^{k-j} \cdot \frac{-\gamma \big( k-j, (-\lambda+\mu)x \big)}{\Gamma(k-j)}\\
    &-&\frac{1}{2} \mu e^{-\mu(x+2y)} \left( \frac{\lambda}{\mu+\lambda} \right)^k \cdot \frac{j}{\lambda+\mu} \cdot \left(1- \frac{\gamma \big( k-j, (\lambda+\mu)x \big)}{\Gamma(k-j)} \right)\\
    &-& \frac{1}{2} \mu e^{-\mu(3x+2y)} \left( \frac{\lambda}{\mu+\lambda} \right)^j \frac{\lambda^{(k-j)}}{\Gamma(k-j) } \\
&\times& \left(  \sum_{i=0}^{k-j-1} (-1)^{i} \cdot \frac{e^{(\mu-\lambda)x}(k-j-1)!}{(k-j-1-i)!\cdot \mu-\lambda)^{i+1}} \cdot x^{k-j-1-i} - (-1)^{k-j-1}\frac{(k-j-1)!}{\mu-\lambda)^{k-j-1}} \right)\\
    &-&\frac{1}{2} \mu e^{-\mu(x+2y)} \left( \frac{\lambda}{\mu+\lambda} \right)^k \cdot \frac{k-j}{\lambda+\mu} \cdot \left(1- \frac{\gamma \big( k-j+1, (\lambda-\mu)x \big)}{\Gamma(k-j+1)} \right).
            \end{eqnarray*}
    \end{enumerate}

    \item when $n<m<n+j < m+k$:
    \begin{eqnarray*}
    \lefteqn{ \mathbb{P} \left( O^{(1)}_{n,n+j} > x, O^{(2)}_{m,m+k} > y  \right) } \\
    &=&\frac{1}{2} \mu e^{-\mu(2x+2y)} \left( \frac{\lambda}{2\mu + \lambda} \right)^{n+j-m} \left( \frac{\lambda}{\mu + \lambda} \right)^{2m+k-2n-j}\\
    &\times& \left( 1 + \frac{n+j-m}{2\mu + \lambda} + \frac{m+k-n-j}{\mu + \lambda} + \frac{1}{2 \mu} \right)
    \end{eqnarray*}
    
    \item when $n<n+j<m < m+k$:
    \begin{eqnarray*}
    \mathbb{P} \left( O^{(1)}_{n,n+j} > x, O^{(2)}_{m,m+k} > y  \right)  
    &=& \frac{1}{4} e^{-\mu(2x+2y)} \left( \frac{\lambda}{\lambda+\mu} \right)^{j+k} \cdot \left( 2 + \frac{2\mu k}{\lambda+\mu} \right)
    \end{eqnarray*}
    
    \item when $n<m < m+k<n+j$:
    \begin{eqnarray*}
    \lefteqn { \mathbb{P} \left( O^{(1)}_{n,n+j} > x, O^{(2)}_{m,m+k} > y  \right) } \\
    &=& \frac{1}{4} e^{-\mu(2x+2y)} \left( \frac{\lambda}{2\mu + \lambda} \right)^{k} \left( \frac{\lambda}{\mu + \lambda} \right)^{j-k} \cdot \left( 2+ \frac{2\mu k}{2\mu + \lambda} \right)
    \end{eqnarray*}
    
    \item when $n = m$, $j > k$:
    \begin{eqnarray*}
   \lefteqn {\mathbb{P} \left( O^{(1)}_{n,n+j} > x, O^{(2)}_{m,m+k} > y  \right) }\\
    &=& \frac{1}{4} e^{-\mu(3x+2y)}  \left( \frac{\lambda}{2\mu + \lambda} \right)^{j-k} \left( \frac{\lambda}{\mu + \lambda} \right)^{k} \cdot \left( 2\mu x +2 + 2\mu \cdot \frac{j-k}{2\mu+\lambda} \right)
    \end{eqnarray*}
    \item when $m < n < m+j < n+k$.
    \begin{eqnarray*}
    \lefteqn{ \mathbb{P} \left( O^{(1)}_{n,n+j} > x, O^{(2)}_{m,m+k} > y  \right) } \\
    &=&\frac{1}{4} e^{-\mu(2x+2y)} \left( \frac{\lambda}{\mu + \lambda} \right)^{2n+j-2m-k} \left( \frac{\lambda}{2\mu + \lambda} \right)^{m-n+k} \cdot \left( 2+ 2\mu\left( \frac{n-m}{\lambda+\mu} + \frac{m-n+k}{2\mu + \lambda} \right) \right)
    \end{eqnarray*}
    \item when $m < m+j < n < n+k$.
    \begin{eqnarray*}
    \lefteqn{ \mathbb{P} \left( O^{(1)}_{n,n+j} > x, O^{(2)}_{m,m+k} > y  \right) } \\
    &=&\frac{1}{4} e^{-\mu(2x+2y)} \left( \frac{\lambda}{\mu + \lambda} \right)^{j+k} \cdot \left( 2+ \frac{ 2\mu k }{\mu + \lambda} \right)
    \end{eqnarray*}
    \item when $m < n < n+k < m+j$.
    \begin{eqnarray*}
    \lefteqn{ \mathbb{P} \left( O^{(1)}_{n,n+j} > x, O^{(2)}_{m,m+k} > y  \right) } \\
    &=& \frac{1}{4} e^{-\mu(2x+2y)} \left( \frac{\lambda}{\mu + \lambda} \right)^{k-j} \left( \frac{\lambda}{2\mu + \lambda} \right)^{j} \cdot \left( 2+ 2\mu\left( \frac{k-j}{\lambda+\mu} + \frac{j}{2\mu + \lambda} \right) \right).
    \end{eqnarray*}
\end{enumerate}
\begin{proof}
 The proof is given in the Appendix.  
\end{proof}
\end{theorem}

\begin{remark}
It is important to note that the other joint probabilities are omitted from the paper, however, they can be easily calculated via the tail distributions of Theorem \ref{ovlp_diff_customers_prop} and the marginal distributions from Theorem \ref{marginals}. 
\end{remark}

\section{Numerical Experiments} \label{numerics}
In this section, we conduct numerical experiments to confirm the theoretical results we obtained in the above sections. Our numerical experiments reflect the different relationships between $\lambda$ and $\mu$, hence we provide three sets of parameters corresponding to the cases $\lambda > \mu$, $\lambda = \mu$, $\lambda < \mu$. \\
\subsection{ $\lambda > \mu$}
For our first numerical experiment, the parameters we use are: 
$$\lambda = 10, \quad \mu = 2, \quad n = 1000000, \quad x = 0.004, \quad y = 0.006.$$
Using the above parameters we obtain the results in Table ~\ref{T:1} below.
\begin{table}[H] 
\begin{center}
\begin{tabular}{|l|l|l|l|l|l|l|} 
\hline
\textbf{}                 & \textbf{Theoretical} & \textbf{Simulated} & \textbf{$j$} & \textbf{$k$}& \textbf{$m-n$} \\ \hline
\textbf{$n = m$, $j = k$} & 0.333052      & 0.33360& 2   & 2   & 0        \\ \hline
\textbf{$n=m, j < k$}              & 0.2688             & 0.2695  & 2    & 6   & 0    \\ \hline
\textbf{$n<m<n+j < m+k$}    & 0.2109             & 0.2110   & 3    & 5   & 1   \\ \hline
\textbf{$n<n+j<m < m+k$}    & 0.2460              & 0.2458  & 2    & 5   & 3    \\ \hline
\textbf{$n<m < m+k<n+j$}    & 0.1240              & 0.1239   & 6    & 2   & 1   \\ \hline
\textbf{$n = m$, $j > k$}    & 0.1362              & 0.1361   & 2    & 6   & 0   \\ \hline
\textbf{$m < n < m+j < n+k$}    & 0.1714              & 0.1717  & 5    & 3   & -1    \\ \hline
\textbf{$m < m+j < n < n+k$}    & 0.2574              & 0.2574  & 3    & 2   & -3    \\ \hline
\textbf{$m < n < n+k < m+j$}    & 0.2305              & 0.2307   & 2   & 6   & -1    \\ \hline
\end{tabular}
\end{center}
\caption{Numerical Experiments for $\lambda > \mu$.}
\label{T:1}
\end{table}

\subsection{ $\lambda = \mu$}
For our second numerical experiment, the parameters we use are:
$$ \lambda = 7, \quad \mu = 7, \quad n = 1000000, \quad x = 0.04, \quad y = 0.06.$$
Using the above parameters we obtain the results in Table \ref{T2} below.
\begin{table}[H]
\begin{tabular}{|l|l|l|l|l|l|l|}
\hline
\textbf{}                 & \textbf{Theoretical} & \textbf{Simulated} & \textbf{$j$} & \textbf{$k$}& \textbf{$m-n$} \\ \hline
\textbf{$n = m$, $j = k$} & 0.03082      & 0.03047  & 2   & 2   & 0        \\ \hline
\textbf{$n=m, j < k$}              & 0.00692           & 0.00691  & 2    & 6   & 0    \\ \hline
\textbf{$n<m<n+j < m+k$}    & 0.00271             & 0.00284   & 3    & 5   & 1   \\ \hline
\textbf{$n<n+j<m < m+k$}    & 0.00373              & 0.00383  & 2    & 5   & 3    \\ \hline
\textbf{$n<m < m+k<n+j$}    & 0.00142              & 0.00133   & 6    & 2   & 1   \\ \hline
\textbf{$n = m$, $j > k$}    & 0.00075              & 0.00077   & 2    & 6   & 0   \\ \hline
\textbf{$m < n < m+j < n+k$}    & 0.00185              & 0.00197  & 5    & 3   & -1    \\ \hline
\textbf{$m < m+j < n < n+k$}    & 0.00770              & 0.00764  & 3    & 2   & -3    \\ \hline
\textbf{$m < n < n+k < m+j$}    & 0.00314              & 0.00332   & 2   & 6   & -1    \\ \hline
\end{tabular} 
\caption{Numerical Experiments for $\lambda = \mu$.}
\label{T2}
\end{table}

\subsection{ $\lambda < \mu$}
For our final numerical experiment, the parameters we use are:
$$ \lambda = 4, \quad \mu = 5, \quad n = 1000000, \quad x = 0.004, \quad y = 0.006.$$
Using the above parameters we obtain the results in Table \ref{T3} below.
\begin{table}[ht]
\begin{tabular}{|l|l|l|l|l|l|l|}
\hline
\textbf{}                 & \textbf{Theoretical} & \textbf{Simulated} & \textbf{$j$} & \textbf{$k$}& \textbf{$m-n$} \\ \hline
\textbf{$n = m$, $j = k$} & 0.08937      & 0.08967  & 2   & 2   & 0        \\ \hline
\textbf{$n=m, j < k$}              & 0.01139           & 0.01140  & 2    & 6   & 0    \\ \hline
\textbf{$n<m<n+j < m+k$}    & 0.00487             & 0.00476   & 3    & 5   & 1   \\ \hline
\textbf{$n<n+j<m < m+k$}    & 0.00585              & 0.00602  & 2    & 5   & 3    \\ \hline
\textbf{$n<m < m+k<n+j$}    & 0.00171              & 0.00151   & 6    & 2   & 1   \\ \hline
\textbf{$n = m$, $j > k$}    & 0.00143              & 0.00125   & 2    & 6   & 0   \\ \hline
\textbf{$m < n < m+j < n+k$}    & 0.00327              & 0.00335  & 5    & 3   & -1    \\ \hline
\textbf{$m < m+j < n < n+k$}    & 0.01756              & 0.01714  & 3    & 2   & -3    \\ \hline
\textbf{$m < n < n+k < m+j$}    & 0.00567              & 0.00558   & 2   & 6   & -1    \\ \hline
\end{tabular} 
\caption{Numerical Experiments for $\lambda < \mu$.}
\label{T3}
\end{table}

We observe in Tables \ref{T:1}, \ref{T2} and \ref{T3} that our simulation experiments match the theoretical results obtained in the paper. This further provides evidence for the accuracy of our theoretical findings.

\section{Conclusion} \label{conclusion}

In this paper, we consider the overlap times for customers in the two dimensional $(M/M/\infty)^2$ tandem queue where the service times of both stations are identically distributed.   We derive explicit formula for the overlap of pairs of customers in both stations.  This analysis is useful for understanding how much a pair of customers overlap in both stations or how much different pairs overlap in each of the stations.   Despite our analysis, there are many avenues for additional research.  First, it is important to generalize our work to an arbitrary number of stations and where the service times at each station are not identically distributed.  This is non-trivial from an integration perspective since for a system with $N$ stations would need to solve a $2N+1$ dimensional integral.  However, we provide this conjecture for the $N$-dimensional case for the joint distribution of overlap times when the service times are identical below.  

\begin{conjecture}
Let $O^{(i)}_{n,n+k}$ be the overlap times for the $n^{th}$ and $(n+k)^{th}$ customers in the $(M/M/\infty)^j$ tandem queue with identically distributed stations, then the joint tail distribution of the overlap times is equal to
\begin{eqnarray*}
\mathbb{P} \left( O^{(1)}_{n,n+k} > x_1, O^{(2)}_{n,n+k} > x_2,  \dots ,  O^{(j)}_{n,n+k} > x_j  \right) &=& \frac{e^{-2\mu \left(\sum^j_{i=1} x_i \right)} }{2^{j-1} } \left( \frac{\lambda}{\lambda+\mu} \right)^{k}  .
\end{eqnarray*}
\end{conjecture}

A second extension would consider generalizing the model to have general arrival and service processes.  In particular, considering generalizations such as dependent arrivals or service times like in \citet{pang2012impact, pang2013two, daw2018queues} would be of great interest, especially in the quantum computing application. In particular, dependent service times would impact the tail distribution quite significantly since it no longer splits into a product of individual distributions.  It would also be interesting to consider the situation where there are a finite number of servers or customer abandonments to assess the impact on the tail distribution of the overlap time, see for example \citet{ko2022overlapping}.  Finally, we are also interested in exploring overlap times in queueing systems with time varying rates.  We hope to consider these important extensions in future work.  



\section*{Acknowledgements}
Jamol Pender would like to acknowledge the gracious support of the National Science Foundation DMS Award \# 2206286.  Valarie Gao would like to acknowledge the gracious support of the ORACL program for supporting her graduate studies.

\appendix
\section{Appendix} \label{appendix}
In this section, we provide the complete proof of Theorem \ref{ovlp_diff_customers}.  We break the proof to each of the nine cases, which represent the ordering of the customers whose overlap we compare.  

\subsection{Proof for CASE I in \ref{ovlp_diff_customers_prop}}
\begin{proof}
    The probability of concern is $\mathbb{P} \left( O^{(1)}_{n,n+j} > x, O^{(2)}_{n,n+k} > y  \right)$ where we assume $j < k$. We separate the calculation into two parts, one with $\mathcal{E}_{(2)}$ and one with $\mathcal{E}_{(2)}^C$.
\end{proof}

\subsection{Proof for CASE II in \ref{ovlp_diff_customers_prop}}
$\mathbb{P} \left( O^{(1)}_{n,n+j} > x, O^{(2)}_{n,n+k} > y  \right)$ with $j < k$.
\subsection{Condition 1:}
\subsubsection{$\lambda > \mu$}
\begin{eqnarray*}
\lefteqn{ \mathbb{P} \left( S^{(1)}_n - \mathcal{A}_{j} > x , S^{(1)}_{n+j} > x, S^{(2)}_{n+k} > y, S^{(1)}_{n} + S^{(2)}_{n} - \mathcal{A}_{j} - \Tilde{\mathcal{A}}_{k-j} - S^{(1)}_{n+k} > y , \mathcal{E}_{(2)}  \right) } \\  
&=& \int_{\mathbb{R}_+^7} \{z_3 > x + z_1 \} \cdot \{ z_5 > x \} \cdot \{ z_7 > y \} \cdot \{ z_3+z_4 - z_1-z_2-z_6>y \} \cdot \{ z_1+z_2+z_6>z_3 \}\\
&\times& \mu^5 e^{-\mu(z_3+z_4+z_5+z_6+z_7)} \cdot \frac{\lambda^{j(k-j)} z_1^{j-1} z_2^{k-j-1}}{\Gamma(k-j) \Gamma(j)} \cdot e^{-\lambda (z_1+z_2)} dz_1 dz_2 dz_3 dz_4 dz_5 dz_6 dz_7 \\
&=& e^{-\mu (x+y)} \int_{\mathbb{R}_+^5} \{z_3 > x + z_1 \} \cdot  \{ z_3+z_4 - z_1-z_2-z_6>y \} \cdot \{ z_1+z_2+z_6>z_3 \}\\
&\times& \mu^3 e^{-\mu(z_3+z_4+z_6)} \cdot  \frac{\lambda^{j(k-j)} z_1^{j-1} z_2^{k-j-1}}{\Gamma(k-j) \Gamma(j)}  \cdot e^{-\lambda (z_1+z_2)} dz_1 dz_2 dz_3 dz_4 dz_6 \\
&=& e^{-\mu (x+2y)} \int_{\mathbb{R}_+^4} \{z_3 > x + z_1 \} \cdot \{ z_1+z_2+z_6>z_3 \}\\
&\times& \mu^2 e^{-\mu(2 z_6+z_2+z_1)} \cdot  \frac{\lambda^{j(k-j)} z_1^{j-1} z_2^{k-j-1}}{\Gamma(k-j) \Gamma(j)}  \cdot e^{-\lambda (z_1+z_2)}  dz_1 dz_2 dz_3 dz_6 \\
&=& \frac{1}{2} e^{-\mu (x+2y)} \int_{\mathbb{R}_+^3} \{z_3 > x + z_1 \} \{z_3 > z_2 + z_1 \} \\
&\times& \mu e^{-\mu(-z_2-z_1+2z_3)} \cdot  \frac{\lambda^{j(k-j)} z_1^{j-1} z_2^{k-j-1}}{\Gamma(k-j) \Gamma(j)}  \cdot e^{-\lambda (z_1+z_2)}  dz_1 dz_2 dz_3 \\
&+&e^{-\mu (x+2y)} \int_{\mathbb{R}_+^3} \{z_3 > x + z_1 \} \cdot \{ z_3< z_1+z_2 \} \\
&\times& \frac{1}{2} \mu e^{-\mu(z_2+z_1)} \cdot  \frac{\lambda^{j(k-j)} z_1^{j-1} z_2^{k-j-1}}{\Gamma(k-j) \Gamma(j)}  \cdot e^{-\lambda (z_1+z_2)}  dz_1 dz_2 dz_3\\
&=& \frac{1}{2} e^{-\mu (x+2y)} \int_{\mathbb{R}_+^3} \{z_3 > x + z_1 \} \{z_2 <x \} \\
&\times& \mu e^{-\mu(-z_2-z_1+2z_3)} \cdot  \frac{\lambda^{j(k-j)} z_1^{j-1} z_2^{k-j-1}}{\Gamma(k-j) \Gamma(j)}  \cdot e^{-\lambda (z_1+z_2)}  dz_1 dz_2 dz_3  \\  
&+& \frac{1}{2} e^{-\mu (x+2y)} \int_{\mathbb{R}_+^3} \{z_3 > z_2 + z_1 \} \{z_2 >x \} \\
&\times&\mu e^{-\mu(-z_2-z_1+2z_3)} \cdot  \frac{\lambda^{j(k-j)} z_1^{j-1} z_2^{k-j-1}}{\Gamma(k-j) \Gamma(j)}  \cdot e^{-\lambda (z_1+z_2)}  dz_1 dz_2 dz_3  \\  
&+&  \frac{1}{2} e^{-\mu (x+2y)} \int_{\mathbb{R}_+^2}\{z_2>x\} (z_2-x) \mu e^{-\mu(z_2+z_1)} \cdot  \frac{\lambda^{j(k-j)} z_1^{j-1} z_2^{k-j-1}}{\Gamma(k-j) \Gamma(j)}  \cdot e^{-\lambda (z_1+z_2)}  dz_1 dz_2 \\
&=& \frac{1}{4} e^{-\mu (3x+2y)} \int_{\mathbb{R}_+^2} \{z_2 <x \}  e^{-(\mu+\lambda)z_1} \cdot  \frac{\lambda^{j(k-j)} z_1^{j-1} z_2^{k-j-1}}{\Gamma(k-j) \Gamma(j)}  \cdot e^{-(\lambda-\mu) z_2}  dz_1 dz_2 \\
&+& \frac{1}{4} e^{-\mu (x+2y)} \int_{\mathbb{R}_+^2} \{z_2 > x \}  e^{-(\mu+\lambda)z_1} \cdot  \frac{\lambda^{j(k-j)} z_1^{j-1} z_2^{k-j-1}}{\Gamma(k-j) \Gamma(j)}  \cdot e^{-(\mu+\lambda)z_2} dz_1 dz_2 \\
&+& \frac{\mu}{2} e^{-\mu (x+2y)} \cdot \frac{\lambda^{j}}{(\mu+\lambda)^{j}} \cdot \frac{\lambda^{k-j}}{(\mu+\lambda)^{k-j}}\left(\frac{k-j}{\mu+\lambda}-x\right)\\
&-&\frac{\mu}{2} e^{-\mu (x+2y)} \cdot \frac{\lambda^{j}}{(\mu+\lambda)^{j}} \cdot \frac{\lambda^{k-j}}{(\mu+\lambda)^{k-j}} \left(\frac{\gamma\left(k-j+1,x(\lambda+\mu)\right)}{(\lambda+\mu)\Gamma(k-j)} -x \frac{\gamma\left(k-j,x(\lambda+\mu)\right)}{\Gamma(k-j)} \right)\\
&=& \frac{1}{4} e^{-\mu (3x+2y)}\cdot \frac{\lambda^{j}}{(\mu+\lambda)^{j}} \int_{0}^x \frac{\lambda^{(k-j)}  z_2^{k-j-1}}{\Gamma(k-j) }  \cdot e^{-(\lambda-\mu) z_2}  dz_2 \\
&+& \frac{1}{4} e^{-\mu (x+2y)}\cdot \frac{\lambda^{j}}{(\mu+\lambda)^{j}} \int_{x}^\infty \frac{\lambda^{(k-j)}  z_2^{k-j-1}}{\Gamma(k-j) }  \cdot e^{-(\lambda+\mu) z_2}  dz_2 \\
&+& \frac{\mu}{2} e^{-\mu (x+2y)} \cdot \frac{\lambda^{j}}{(\mu+\lambda)^{j}} \cdot \frac{\lambda^{k-j}}{(\mu+\lambda)^{k-j}}\left(\frac{k-j}{\mu+\lambda}-x\right)\\
&-&\frac{\mu}{2} e^{-\mu (x+2y)} \cdot \frac{\lambda^{j}}{(\mu+\lambda)^{j}} \cdot \frac{\lambda^{k-j}}{(\mu+\lambda)^{k-j}} \left(\frac{\gamma\left(k-j+1,x(\lambda+\mu)\right)}{(\lambda+\mu)\Gamma(k-j)} -x \frac{\gamma\left(k-j,x(\lambda+\mu)\right)}{\Gamma(k-j)} \right)\\
&=& \frac{1}{4} e^{-\mu (3x+2y)}\cdot \frac{\lambda^{j}}{(\mu+\lambda)^{j}} \cdot \frac{\lambda^{k-j}}{(\lambda - \mu)^{k-j}} \cdot \frac{1}{\Gamma(k-j)} \cdot \gamma (k-j,(\lambda-\mu)x)\\
&+& \frac{1}{4} e^{-\mu (x+2y)}\cdot \frac{\lambda^{j}}{(\mu+\lambda)^{j}} \cdot \frac{\lambda^{k-j}}{(\lambda +\mu)^{k-j}} \cdot \left(1- \frac{1}{\Gamma(k-j)} \cdot \gamma (k-j,(\lambda+\mu)x)\right) \\
&+& \frac{\mu}{2} e^{-\mu (x+2y)} \cdot \frac{\lambda^{j}}{(\mu+\lambda)^{j}} \cdot \frac{\lambda^{k-j}}{(\mu+\lambda)^{k-j}}\left(\frac{k-j}{\mu+\lambda}-x\right)\\
&-&\frac{\mu}{2} e^{-\mu (x+2y)} \cdot \frac{\lambda^{j}}{(\mu+\lambda)^{j}} \cdot \frac{\lambda^{k-j}}{(\mu+\lambda)^{k-j}} \left(\frac{\gamma\left(k-j+1,x(\lambda+\mu)\right)}{(\lambda+\mu)\Gamma(k-j)} -x \frac{\gamma\left(k-j,x(\lambda+\mu)\right)}{\Gamma(k-j)} \right)
\end{eqnarray*}

\subsubsection{$\lambda = \mu$}
\begin{eqnarray*}
\lefteqn{ \mathbb{P} \left( S^{(1)}_n - \mathcal{A}_{j} > x , S^{(1)}_{n+j} > x, S^{(2)}_{n+k} > y, S^{(1)}_{n} + S^{(2)}_{n} - \mathcal{A}_{j} - \Tilde{\mathcal{A}}_{k-j} - S^{(1)}_{n+k} > y , \mathcal{E}_{(2)}  \right) } \\  
&=& \int_{\mathbb{R}_+^7} \{z_3 > x + z_1 \} \cdot \{ z_5 > x \} \cdot \{ z_7 > y \} \cdot \{ z_3+z_4 - z_1-z_2-z_6>y \} \cdot \{ z_1+z_2+z_6>z_3 \}\\
&\times& \mu^5 e^{-\mu(z_3+z_4+z_5+z_6+z_7)} \cdot \frac{\lambda^{j(k-j)} z_1^{j-1} z_2^{k-j-1}}{\Gamma(k-j) \Gamma(j)} \cdot e^{-\lambda (z_1+z_2)} dz_1 dz_2 dz_3 dz_4 dz_5 dz_6 dz_7 \\
&=& e^{-\mu (x+y)} \int_{\mathbb{R}_+^5} \{z_3 > x + z_1 \} \cdot  \{ z_3+z_4 - z_1-z_2-z_6>y \} \cdot \{ z_1+z_2+z_6>z_3 \}\\
&\times& \mu^3 e^{-\mu(z_3+z_4+z_6)} \cdot  \frac{\lambda^{j(k-j)} z_1^{j-1} z_2^{k-j-1}}{\Gamma(k-j) \Gamma(j)}  \cdot e^{-\lambda (z_1+z_2)} dz_1 dz_2 dz_3 dz_4 dz_6 \\
&=& e^{-\mu (x+2y)} \int_{\mathbb{R}_+^4} \{z_3 > x + z_1 \} \cdot \{ z_1+z_2+z_6>z_3 \}\\
&\times& \mu^2 e^{-\mu(2 z_6+z_2+z_1)} \cdot  \frac{\lambda^{j(k-j)} z_1^{j-1} z_2^{k-j-1}}{\Gamma(k-j) \Gamma(j)}  \cdot e^{-\lambda (z_1+z_2)}  dz_1 dz_2 dz_3 dz_6 \\
&=& \frac{1}{2} e^{-\mu (x+2y)} \int_{\mathbb{R}_+^3} \{z_3 > x + z_1 \} \{z_3 > z_2 + z_1 \} \\
&\times& \mu e^{-\mu(-z_2-z_1+2z_3)} \cdot  \frac{\lambda^{j(k-j)} z_1^{j-1} z_2^{k-j-1}}{\Gamma(k-j) \Gamma(j)}  \cdot e^{-\lambda (z_1+z_2)}  dz_1 dz_2 dz_3 \\
&+&e^{-\mu (x+2y)} \int_{\mathbb{R}_+^3} \{z_3 > x + z_1 \} \cdot \{ z_3< z_1+z_2 \} \\
&\times& \frac{1}{2} \mu e^{-\mu(z_2+z_1)} \cdot  \frac{\lambda^{j(k-j)} z_1^{j-1} z_2^{k-j-1}}{\Gamma(k-j) \Gamma(j)}  \cdot e^{-\lambda (z_1+z_2)}  dz_1 dz_2 dz_3\\
&=& \frac{1}{2} e^{-\mu (x+2y)} \int_{\mathbb{R}_+^3} \{z_3 > x + z_1 \} \{z_2 <x \} \\
&\times& \mu e^{-\mu(-z_2-z_1+2z_3)} \cdot  \frac{\lambda^{j(k-j)} z_1^{j-1} z_2^{k-j-1}}{\Gamma(k-j) \Gamma(j)}  \cdot e^{-\lambda (z_1+z_2)}  dz_1 dz_2 dz_3  \\  
&+& \frac{1}{2} e^{-\mu (x+2y)} \int_{\mathbb{R}_+^3} \{z_3 > z_2 + z_1 \} \{z_2 >x \} \\
&\times& \mu e^{-\mu(-z_2-z_1+2z_3)} \cdot  \frac{\lambda^{j(k-j)} z_1^{j-1} z_2^{k-j-1}}{\Gamma(k-j) \Gamma(j)}  \cdot e^{-\lambda (z_1+z_2)}  dz_1 dz_2 dz_3  \\  
&+&  \frac{1}{2} e^{-\mu (x+2y)} \int_{\mathbb{R}_+^2}\{z_2>x\} (z_2-x) \mu e^{-\mu(z_2+z_1)} \cdot  \frac{\lambda^{j(k-j)} z_1^{j-1} z_2^{k-j-1}}{\Gamma(k-j) \Gamma(j)}  \cdot e^{-\lambda (z_1+z_2)}  dz_1 dz_2 \\
&=& \frac{1}{4} e^{-\mu (3x+2y)} \int_{\mathbb{R}_+^2} \{z_2 <x \}  e^{-(\mu+\lambda)z_1} \cdot  \frac{\lambda^{j(k-j)} z_1^{j-1} z_2^{k-j-1}}{\Gamma(k-j) \Gamma(j)}  dz_1 dz_2 \\
&+& \frac{1}{4} e^{-\mu (x+2y)} \int_{\mathbb{R}_+^2} \{z_2 > x \}  e^{-(\mu+\lambda)z_1} \cdot  \frac{\lambda^{j(k-j)} z_1^{j-1} z_2^{k-j-1}}{\Gamma(k-j) \Gamma(j)}  \cdot e^{-(\mu+\lambda)z_2} dz_1 dz_2 \\
&+& \frac{\mu}{2} e^{-\mu (x+2y)} \cdot \frac{\lambda^{j}}{(\mu+\lambda)^{j}} \cdot \frac{\lambda^{k-j}}{(\mu+\lambda)^{k-j}}\left(\frac{k-j}{\mu+\lambda}-x\right)\\
&-&\frac{\mu}{2} e^{-\mu (x+2y)} \cdot \frac{\lambda^{j}}{(\mu+\lambda)^{j}} \cdot \frac{\lambda^{k-j}}{(\mu+\lambda)^{k-j}} \left(\frac{\gamma\left(k-j+1,x(\lambda+\mu)\right)}{(\lambda+\mu)\Gamma(k-j)} -x \frac{\gamma\left(k-j,x(\lambda+\mu)\right)}{\Gamma(k-j)} \right)\\
&=& \frac{1}{4} e^{-\mu (3x+2y)}\cdot \frac{\lambda^{j}}{(\mu+\lambda)^{j}} \int_{0}^x \frac{\lambda^{(k-j)}  z_2^{k-j-1}}{\Gamma(k-j) }   dz_2 \\
&+& \frac{1}{4} e^{-\mu (x+2y)}\cdot \frac{\lambda^{j}}{(\mu+\lambda)^{j}} \int_{x}^\infty \frac{\lambda^{(k-j)}  z_2^{k-j-1}}{\Gamma(k-j) }  \cdot e^{-(\lambda+\mu) z_2}  dz_2 \\
&+& \frac{\mu}{2} e^{-\mu (x+2y)} \cdot \frac{\lambda^{j}}{(\mu+\lambda)^{j}} \cdot \frac{\lambda^{k-j}}{(\mu+\lambda)^{k-j}}\left(\frac{k-j}{\mu+\lambda}-x\right)\\
&-&\frac{\mu}{2} e^{-\mu (x+2y)} \cdot \frac{\lambda^{j}}{(\mu+\lambda)^{j}} \cdot \frac{\lambda^{k-j}}{(\mu+\lambda)^{k-j}} \left(\frac{\gamma\left(k-j+1,x(\lambda+\mu)\right)}{(\lambda+\mu)\Gamma(k-j)} -x \frac{\gamma\left(k-j,x(\lambda+\mu)\right)}{\Gamma(k-j)} \right)\\
&=& \frac{1}{4} e^{-\mu (3x+2y)}\cdot \frac{\lambda^{j}}{(\mu+\lambda)^{j}} \cdot \frac{x^{k-j}}{{k-j}} \cdot \frac{\lambda^{(k-j)}}{\Gamma(k-j) }  \\
&+& \frac{1}{4} e^{-\mu (x+2y)}\cdot \frac{\lambda^{j}}{(\mu+\lambda)^{j}} \cdot \frac{\lambda^{k-j}}{(\lambda +\mu)^{k-j}} \cdot \left(1- \frac{1}{\Gamma(k-j)} \cdot \gamma (k-j,(\lambda+\mu)x)\right) \\
&+& \frac{\mu}{2} e^{-\mu (x+2y)} \cdot \frac{\lambda^{j}}{(\mu+\lambda)^{j}} \cdot \frac{\lambda^{k-j}}{(\mu+\lambda)^{k-j}}\left(\frac{k-j}{\mu+\lambda}-x\right)\\
&-&\frac{\mu}{2} e^{-\mu (x+2y)} \cdot \frac{\lambda^{j}}{(\mu+\lambda)^{j}} \cdot \frac{\lambda^{k-j}}{(\mu+\lambda)^{k-j}} \left(\frac{\gamma\left(k-j+1,x(\lambda+\mu)\right)}{(\lambda+\mu)\Gamma(k-j)} -x \frac{\gamma\left(k-j,x(\lambda+\mu)\right)}{\Gamma(k-j)} \right)
\end{eqnarray*}

\subsubsection{$\lambda < \mu$} \label{A33}

\begin{lemma}
    For $a>0$, $s>0$ and $t > 0$, we have:
    \[ \int_0^x t^{s-1}e^{at} dt = \sum_{i=0}^{s-1} (-1)^{i} \cdot \frac{e^{ax}(s-1)!}{(s-1-i)!\cdot a^{i+1}} \cdot x^{s-1-i} - (-1)^{s-1}\frac{(s-1)!}{a^{s-1}} \]
\end{lemma}
\begin{proof}
    Define $f(t) = t^{s-1}$ and $dg = e^{at} dt$, then we have:
    \begin{align} \nonumber
        f &= t^{s-1} &dg &=e^{at} dt\\ \nonumber
        df &= (s-1)t^{s-2} dt &g &= \frac{e^{at}}{a}
    \end{align}
    By Integral by parts we have:
    \begin{align} \nonumber
        \int t^{s-1}e^{at} dt &= t^{s-1} \cdot \frac{e^{at}}{a} - \frac{s-1}{a} \int  t^{s-2}e^{at} dt
    \end{align}
    Notice that by reapplying integral by parts to $\int  t^{s-i}e^{at} dt$ $s$ times we can get:
    \begin{align} \nonumber
        \int t^{s-1}e^{at} dt &= t^{s-1} \cdot \frac{e^{at}}{a} - \frac{s-1}{a}\left( t^{s-2} \cdot \frac{e^{at}}{a} \left(... - \frac{2}{a} (t \cdot \frac{e^{at}}{a} - \frac{1}{a} \frac{e^{at}}{a}) \right) \right) dt\\ \nonumber
        &= t^{s-1} \cdot \frac{e^{at}}{a} - \frac{s-1}{a} \cdot  t^{s-2}  \cdot \frac{e^{at}}{a} + ... + (-1)^{s-1} \cdot \frac{(s-1)!}{a^{s-1}} \cdot \frac{e^{at}}{a}\\ \nonumber
        &= \sum_{i=0}^{s-1} (-1)^{i} \cdot \frac{e^{at}(s-1)!}{(s-1-i)!\cdot a^{i+1}} \cdot t^{s-1-i} - (-1)^{s-1}\frac{(s-1)!}{a^{s-1}}
    \end{align}
    Therefore
    \begin{align} \nonumber
        \int_0^x t^{s-1}e^{at} dt &=  \sum_{i=0}^{s-1} (-1)^{i} \cdot \frac{e^{as}(s-1)!}{(s-1-i)!\cdot a^{i+1}} \cdot t^{s-1-i} -(-1)^{s-1}\frac{(s-1)!}{a^{s-1}} 
    \end{align}
\end{proof}

\begin{eqnarray*}
\lefteqn{ \mathbb{P} \left( S^{(1)}_n - \mathcal{A}_{j} > x , S^{(1)}_{n+j} > x, S^{(2)}_{n+k} > y, S^{(1)}_{n} + S^{(2)}_{n} - \mathcal{A}_{j} - \Tilde{\mathcal{A}}_{k-j} - S^{(1)}_{n+k} > y , \mathcal{E}_{(2)}  \right) } \\  
&=& \int_{\mathbb{R}_+^7} \{z_3 > x + z_1 \} \cdot \{ z_5 > x \} \cdot \{ z_7 > y \} \cdot \{ z_3+z_4 - z_1-z_2-z_6>y \} \cdot \{ z_1+z_2+z_6>z_3 \}\\
&\times& \mu^5 e^{-\mu(z_3+z_4+z_5+z_6+z_7)} \cdot \frac{\lambda^{j(k-j)} z_1^{j-1} z_2^{k-j-1}}{\Gamma(k-j) \Gamma(j)} \cdot e^{-\lambda (z_1+z_2)} dz_1 dz_2 dz_3 dz_4 dz_5 dz_6 dz_7 \\
&=& e^{-\mu (x+y)} \int_{\mathbb{R}_+^5} \{z_3 > x + z_1 \} \cdot  \{ z_3+z_4 - z_1-z_2-z_6>y \} \cdot \{ z_1+z_2+z_6>z_3 \}\\
&\times& \mu^3 e^{-\mu(z_3+z_4+z_6)} \cdot  \frac{\lambda^{j(k-j)} z_1^{j-1} z_2^{k-j-1}}{\Gamma(k-j) \Gamma(j)}  \cdot e^{-\lambda (z_1+z_2)} dz_1 dz_2 dz_3 dz_4 dz_6 \\
&=& e^{-\mu (x+2y)} \int_{\mathbb{R}_+^4} \{z_3 > x + z_1 \} \cdot \{ z_1+z_2+z_6>z_3 \}\\
&\times& \mu^2 e^{-\mu(2 z_6+z_2+z_1)} \cdot  \frac{\lambda^{j(k-j)} z_1^{j-1} z_2^{k-j-1}}{\Gamma(k-j) \Gamma(j)}  \cdot e^{-\lambda (z_1+z_2)}  dz_1 dz_2 dz_3 dz_6 \\
&=& \frac{1}{2} e^{-\mu (x+2y)} \int_{\mathbb{R}_+^3} \{z_3 > x + z_1 \} \{z_3 > z_2 + z_1 \} \\
&\times& \mu e^{-\mu(-z_2-z_1+2z_3)} \cdot  \frac{\lambda^{j(k-j)} z_1^{j-1} z_2^{k-j-1}}{\Gamma(k-j) \Gamma(j)}  \cdot e^{-\lambda (z_1+z_2)}  dz_1 dz_2 dz_3 \\
&+&e^{-\mu (x+2y)} \int_{\mathbb{R}_+^3} \{z_3 > x + z_1 \} \cdot \{ z_3< z_1+z_2 \} \\
&\times& \frac{1}{2} \mu e^{-\mu(z_2+z_1)} \cdot  \frac{\lambda^{j(k-j)} z_1^{j-1} z_2^{k-j-1}}{\Gamma(k-j) \Gamma(j)}  \cdot e^{-\lambda (z_1+z_2)}  dz_1 dz_2 dz_3\\
&=& \frac{1}{2} e^{-\mu (x+2y)} \int_{\mathbb{R}_+^3} \{z_3 > x + z_1 \} \{z_2 <x \} \\
&\times& \mu e^{-\mu(-z_2-z_1+2z_3)} \cdot  \frac{\lambda^{j(k-j)} z_1^{j-1} z_2^{k-j-1}}{\Gamma(k-j) \Gamma(j)}  \cdot e^{-\lambda (z_1+z_2)}  dz_1 dz_2 dz_3  \\  
&+& \frac{1}{2} e^{-\mu (x+2y)} \int_{\mathbb{R}_+^3} \{z_3 > z_2 + z_1 \} \{z_2 >x \} \\
&\times& \mu e^{-\mu(-z_2-z_1+2z_3)} \cdot  \frac{\lambda^{j(k-j)} z_1^{j-1} z_2^{k-j-1}}{\Gamma(k-j) \Gamma(j)}  \cdot e^{-\lambda (z_1+z_2)}  dz_1 dz_2 dz_3  \\  
&+&  \frac{1}{2} e^{-\mu (x+2y)} \int_{\mathbb{R}_+^2}\{z_2>x\} (z_2-x) \mu e^{-\mu(z_2+z_1)} \cdot  \frac{\lambda^{j(k-j)} z_1^{j-1} z_2^{k-j-1}}{\Gamma(k-j) \Gamma(j)}  \cdot e^{-\lambda (z_1+z_2)}  dz_1 dz_2 \\
&=& \frac{1}{4} e^{-\mu (3x+2y)} \int_{\mathbb{R}_+^2} \{z_2 <x \}  e^{-(\mu+\lambda)z_1} \cdot  \frac{\lambda^{j(k-j)} z_1^{j-1} z_2^{k-j-1}}{\Gamma(k-j) \Gamma(j)}  \cdot e^{-(\lambda-\mu) z_2}  dz_1 dz_2 \\
&+& \frac{1}{4} e^{-\mu (x+2y)} \int_{\mathbb{R}_+^2} \{z_2 > x \}  e^{-(\mu+\lambda)z_1} \cdot  \frac{\lambda^{j(k-j)} z_1^{j-1} z_2^{k-j-1}}{\Gamma(k-j) \Gamma(j)}  \cdot e^{-(\mu+\lambda)z_2} dz_1 dz_2 \\
&+& \frac{\mu}{2} e^{-\mu (x+2y)} \cdot \frac{\lambda^{j}}{(\mu+\lambda)^{j}} \cdot \frac{\lambda^{k-j}}{(\mu+\lambda)^{k-j}}\left(\frac{k-j}{\mu+\lambda}-x\right)\\
&-&\frac{\mu}{2} e^{-\mu (x+2y)} \cdot \frac{\lambda^{j}}{(\mu+\lambda)^{j}} \cdot \frac{\lambda^{k-j}}{(\mu+\lambda)^{k-j}} \left(\frac{\gamma\left(k-j+1,x(\lambda+\mu)\right)}{(\lambda+\mu)\Gamma(k-j)} -x \frac{\gamma\left(k-j,x(\lambda+\mu)\right)}{\Gamma(k-j)} \right)\\
&=& \frac{1}{4} e^{-\mu (3x+2y)}\cdot \frac{\lambda^{j}}{(\mu+\lambda)^{j}} \int_{0}^x \frac{\lambda^{(k-j)}  z_2^{k-j-1}}{\Gamma(k-j) }  \cdot e^{-(\lambda-\mu) z_2}  dz_2 \\
&+& \frac{1}{4} e^{-\mu (x+2y)}\cdot \frac{\lambda^{j}}{(\mu+\lambda)^{j}} \int_{x}^\infty \frac{\lambda^{(k-j)}  z_2^{k-j-1}}{\Gamma(k-j) }  \cdot e^{-(\lambda+\mu) z_2}  dz_2 \\
&+& \frac{\mu}{2} e^{-\mu (x+2y)} \cdot \frac{\lambda^{j}}{(\mu+\lambda)^{j}} \cdot \frac{\lambda^{k-j}}{(\mu+\lambda)^{k-j}}\left(\frac{k-j}{\mu+\lambda}-x\right)\\
&-&\frac{\mu}{2} e^{-\mu (x+2y)} \cdot \frac{\lambda^{j}}{(\mu+\lambda)^{j}} \cdot \frac{\lambda^{k-j}}{(\mu+\lambda)^{k-j}} \left(\frac{\gamma\left(k-j+1,x(\lambda+\mu)\right)}{(\lambda+\mu)\Gamma(k-j)} -x \frac{\gamma\left(k-j,x(\lambda+\mu)\right)}{\Gamma(k-j)} \right)\\
&=& \frac{1}{4} e^{-\mu (3x+2y)}\cdot \frac{\lambda^{j}}{(\mu+\lambda)^{j}} \frac{\lambda^{(k-j)}}{\Gamma(k-j) } \\
&\times& \left(  \sum_{i=0}^{k-j-1} (-1)^{i} \cdot \frac{e^{(\mu-\lambda)x}(k-j-1)!}{(k-j-1-i)!\cdot \mu-\lambda)^{i+1}} \cdot x^{k-j-1-i} - (-1)^{k-j-1}\frac{(k-j-1)!}{\mu-\lambda)^{k-j-1}} \right)\\
&+& \frac{1}{4} e^{-\mu (x+2y)}\cdot \frac{\lambda^{j}}{(\mu+\lambda)^{j}} \cdot \frac{\lambda^{k-j}}{(\lambda +\mu)^{k-j}} \cdot \left(1- \frac{1}{\Gamma(k-j)} \cdot \gamma (k-j,(\lambda+\mu)x)\right) \\
&+& \frac{\mu}{2} e^{-\mu (x+2y)} \cdot \frac{\lambda^{j}}{(\mu+\lambda)^{j}} \cdot \frac{\lambda^{k-j}}{(\mu+\lambda)^{k-j}}\left(\frac{k-j}{\mu+\lambda}-x\right)\\
&-&\frac{\mu}{2} e^{-\mu (x+2y)} \cdot \frac{\lambda^{j}}{(\mu+\lambda)^{j}} \cdot \frac{\lambda^{k-j}}{(\mu+\lambda)^{k-j}} \left(\frac{\gamma\left(k-j+1,x(\lambda+\mu)\right)}{(\lambda+\mu)\Gamma(k-j)} -x \frac{\gamma\left(k-j,x(\lambda+\mu)\right)}{\Gamma(k-j)} \right)
\end{eqnarray*}

\subsection{Condition 2}
\subsubsection{$\lambda > \mu$} \label{seca41}
\begin{eqnarray*}
    \lefteqn{ \mathbb{P} \left( S^{(1)}_n - \mathcal{A}_{j} > x , S^{(1)}_{n+j} > x, S^{(2)}_{n+k} > y, S^{(1)}_{n} + S^{(2)}_{n} - \mathcal{A}_{j} - \Tilde{\mathcal{A}}_{k-j} - S^{(1)}_{n+k} > y , \mathcal{E}^C_{(2)}  \right) } \\ 
    &=& \int_{\mathbb{R}_+^7} \{z_3 > x + z_1 \} \cdot \{ z_5 > x \} \cdot \{ z_4 > y \} \cdot \{ z_6+z_7 + z_1+z_2-z_3>y \} \cdot \{ z_1+z_2+z_6<z_3 \}\\
    &\times& \mu^5 e^{-\mu(z_3+z_4+z_5+z_6+z_7)} \cdot \frac{\lambda^{j + (k-j)} z_1^{j-1} z_2^{k-j-1}}{\Gamma(k-j) \Gamma(j)} \cdot e^{-\lambda (z_1+z_2)} dz_1 dz_2 dz_3 dz_4 dz_5 dz_6 dz_7 \\
    &=& e^{-\mu (x+y)} \int_{\mathbb{R}_+^5} \{z_3 > x + z_1 \} \cdot  \{ z_6+z_7 + z_1+z_2-z_3>y \} \cdot \{ z_1+z_2+z_6<z_3 \}\\
    &\times& \mu^3 e^{-\mu(z_3+z_6+z_7)} \cdot  \frac{\lambda^{j+(k-j)} z_1^{j-1} z_2^{k-j-1}}{\Gamma(k-j) \Gamma(j)}  \cdot e^{-\lambda (z_1+z_2)} dz_1 dz_2 dz_3 dz_6 dz_7 \\
    &=& e^{-\mu (x+2y)} \int_{\mathbb{R}_+^4} \{z_3 > x + z_1 \} \cdot \{ z_6 < z_3 - z_1-z_2\} \cdot \{ z_3> z_1+z_2\} \cdot \mu^2 \cdot e^{-\mu(2z_3-z_2-z_1)} \\
    &\times& \frac{\lambda^{j+(k-j)} z_1^{j-1} z_2^{k-j-1}}{\Gamma(k-j) \Gamma(j)} \cdot e^{-\lambda (z_1+z_2)} dz_1 dz_2 dz_3 dz_6\\
    &=&e^{-\mu (x+2y)} \int_{\mathbb{R}_+^3} \{z_3 > x + z_1 \} \cdot \{ z_3> z_1+z_2\} \cdot \mu^2 \cdot e^{-\mu(2z_3-z_2-z_1)}\\
    &\times& \frac{\lambda^{j+(k-j)} z_1^{j-1} z_2^{k-j-1}}{\Gamma(k-j) \Gamma(j)} \cdot e^{-\lambda (z_1+z_2)} \cdot (z_3 - z_1 -z_2) dz_1 dz_2 dz_3.
\end{eqnarray*}
To simplify the analysis, we now define $M$ to be $M = \frac{\lambda^{j+(k-j)} z_1^{j-1} z_2^{k-j-1}}{\Gamma(k-j) \Gamma(j)}$. This allows us to split the above equation into three components labeled $A,B$, and $C$ for convenience.
\begin{enumerate}
    \item \begin{eqnarray*}
        \lefteqn{ A = e^{-\mu (x+2y)} \int_{\mathbb{R}_+^3} \{z_3 > x + z_1 \} \cdot \{ z_3> z_1+z_2\} }\\ 
        &\times& \mu^2 \cdot e^{-\mu(2z_3-z_2-z_1)} \cdot M \cdot e^{-\lambda (z_1+z_2)} \cdot z_3 dz_1 dz_2 dz_3\\
        &=& e^{-\mu (x+2y)} \int_{\mathbb{R}_+^3} \{z_3 > x + z_1 \} \cdot \{ x>z_2\} \\
        &\times& \mu^2 \cdot e^{-\mu(2z_3-z_2-z_1)} \cdot M \cdot e^{-\lambda (z_1+z_2)} \cdot z_3 dz_1 dz_2 dz_3\\
        &+& e^{-\mu (x+2y)} \int_{\mathbb{R}_+^3} \{z_3 > z_2 + z_1 \} \cdot \{ x<z_2\} \\
        &\times& \mu^2 \cdot e^{-\mu(2z_3-z_2-z_1)} \cdot M \cdot e^{-\lambda (z_1+z_2)} \cdot z_3 dz_1 dz_2 dz_3\\
        &=& \frac{1}{4} e^{-\mu (3x+2y)} \cdot (2\mu x+1) \int_{\mathbb{R}_+^2} \{ x > z_2\} \cdot e^{-(\mu+\lambda) z_1} \cdot e^{-(\lambda-\mu) z_2} \cdot M dz_1 dz_2 \\
        &+& \frac{1}{4} e^{-\mu (3x+2y)} \cdot (2\mu) \int_{\mathbb{R}_+^2} z_1 \cdot \{ x > z_2\} \cdot e^{-(\mu+\lambda) z_1} \cdot e^{-(\lambda-\mu) z_2} \cdot M dz_1 dz_2 \\
        &+& \frac{1}{4} e^{-\mu (x+2y)}  \int_{\mathbb{R}_+^2} \{ x < z_2\} \cdot e^{-(\mu+\lambda) z_1} \cdot e^{-(\mu+\lambda) z_2} \cdot M dz_1 dz_2 \\
        &+& \frac{1}{4} e^{-\mu (x+2y)} \cdot (2\mu) \int_{\mathbb{R}_+^2} (z_1+z_2) \cdot \{ x < z_2\} \cdot e^{-(\mu+\lambda) z_1} \cdot e^{-(\mu+\lambda) z_2} \cdot M dz_1 dz_2 \\
        &=&  \frac{1}{4} e^{-\mu (3x+2y)} \left( \frac{\lambda}{\mu+\lambda} \right)^j \cdot \left( \frac{\lambda}{\lambda-\mu} \right)^{k-j} \frac{\gamma \big( k-j, (\lambda-\mu)x \big)}{\Gamma(k-j)} \cdot \left( 2\mu x +1 - \frac{2 \mu j}{\mu+\lambda} \right)\\
        &+& \frac{1}{4} e^{-\mu (x+2y)} \cdot \left( \frac{\lambda}{\mu+\lambda} \right)^k \cdot \left( 1- \frac{\gamma \big( k-j, (\lambda+\mu)x \big)}{\Gamma(k-j)} \right) \cdot \left( 1+ \frac{2\mu j}{\mu + \lambda} \right)\\
        &+&\frac{1}{4} e^{-\mu (x+2y)} \cdot \left( \frac{\lambda}{\mu+\lambda} \right)^k \cdot \frac{k-j}{\lambda+\mu} \cdot \left( 1- \frac{\gamma \big( k-j+1, (\lambda+\mu)x \big)}{\Gamma(k-j)} \right). \\
    \end{eqnarray*}
        
    \item \begin{eqnarray*}
        \lefteqn{ B= e^{-\mu (x+2y)} \int_{\mathbb{R}_+^3} \{z_3 > x + z_1 \} \cdot \{ z_3> z_1+z_2\}  }\\ 
        &\times& \mu^2 \cdot e^{-\mu(2z_3-z_2-z_1)} \cdot M \cdot e^{-\lambda (z_1+z_2)} \cdot z_1 dz_1 dz_2 dz_3 \\
        &=& e^{-\mu (x+2y)} \int_{\mathbb{R}_+^3} \{z_3 > x + z_1 \} \cdot \{ x>z_2\} \\
        &\times& \mu^2 \cdot e^{-\mu(2z_3-z_2-z_1)} \cdot M \cdot e^{-\lambda (z_1+z_2)} \cdot z_1 dz_1 dz_2 dz_3\\
        &+& e^{-\mu (x+2y)} \int_{\mathbb{R}_+^3} \{z_3 > z_2 + z_1 \} \cdot \{ x<z_2\} \\
        &\times& \mu^2 \cdot e^{-\mu(2z_3-z_2-z_1)} \cdot M \cdot e^{-\lambda (z_1+z_2)} \cdot z_1 dz_1 dz_2 dz_3\\
        &=& \frac{1}{2} \mu e^{-\mu(3x+2y)}  \int_{\mathbb{R}_+^2} \{x > z_2\} \cdot e^{-(\lambda+\mu)z_1} \cdot e^{-(\lambda-\mu)z_2} \cdot M \cdot z_1 dz_1 dz_2\\
        &+& \frac{1}{2} \mu e^{-\mu(x+2y)}  \int_{\mathbb{R}_+^2} \{x < z_2\} \cdot e^{-(\lambda+\mu)z_1} \cdot e^{-(\lambda+\mu)z_2} \cdot M \cdot z_1 dz_1 dz_2\\
        &=& \frac{1}{2} \mu e^{-\mu(3x+2y)} \left( \frac{\lambda}{\mu+\lambda} \right)^j \cdot \frac{j}{\lambda+\mu} \cdot \left( \frac{\lambda}{\lambda-\mu} \right)^{k-j} \cdot \frac{\gamma \big( k-j, (\lambda-\mu)x \big)}{\Gamma(k-j)}\\
        &+&\frac{1}{2} \mu e^{-\mu(x+2y)} \left( \frac{\lambda}{\mu+\lambda} \right)^k \cdot \frac{j}{\lambda+\mu} \cdot \left(1- \frac{\gamma \big( k-j, (\lambda+\mu)x \big)}{\Gamma(k-j)} \right).\\
    \end{eqnarray*}

    \item \begin{eqnarray*}
        \lefteqn{ C= e^{-\mu (x+2y)} \int_{\mathbb{R}_+^3} \{z_3 > x + z_1 \} \cdot \{ z_3> z_1+z_2\}  }\\ 
        &\times& \mu^2 \cdot e^{-\mu(2z_3-z_2-z_1)} \cdot M \cdot e^{-\lambda (z_1+z_2)} \cdot z_2 dz_1 dz_2 dz_3\\
        &=& e^{-\mu (x+2y)} \int_{\mathbb{R}_+^3} \{z_3 > x + z_1 \} \cdot \{ x>z_2\} \\
        &\times& \mu^2 \cdot e^{-\mu(2z_3-z_2-z_1)} \cdot M \cdot e^{-\lambda (z_1+z_2)} \cdot z_2 dz_1 dz_2 dz_3\\
        &+& e^{-\mu (x+2y)} \int_{\mathbb{R}_+^3} \{z_3 > z_2 + z_1 \} \cdot \{ x<z_2\} \\
        &\times& \mu^2 \cdot e^{-\mu(2z_3-z_2-z_1)} \cdot M \cdot e^{-\lambda (z_1+z_2)} \cdot z_2 dz_1 dz_2 dz_3\\
        &=& \frac{1}{2} \mu e^{-\mu(3x+2y)} \left( \frac{\lambda}{\mu+\lambda} \right)^j \cdot \frac{k-j}{\lambda-\mu} \cdot \left( \frac{\lambda}{\lambda-\mu} \right)^{k-j} \cdot \frac{\gamma \big( k-j+1, (\lambda-\mu)x \big)}{\Gamma(k-j+1)}\\
        &+&\frac{1}{2} \mu e^{-\mu(x+2y)} \left( \frac{\lambda}{\mu+\lambda} \right)^k \cdot \frac{k-j}{\lambda+\mu} \cdot \left(1- \frac{\gamma \big( k-j+1, (\lambda+\mu)x \big)}{\Gamma(k-j+1)} \right).\\
    \end{eqnarray*}
\end{enumerate}
Putting all three parts together we have:
\begin{eqnarray*}
    \lefteqn{ \mathbb{P} \left( S^{(1)}_n - \mathcal{A}_{j} > x , S^{(1)}_{n+j} > x, S^{(2)}_{n+k} > y, S^{(1)}_{n} + S^{(2)}_{n} - \mathcal{A}_{j} - \Tilde{\mathcal{A}}_{k-j} - S^{(1)}_{n+k} > y , \mathcal{E}^C_{(2)}  \right) } \\ 
    &=& A-B-C \\
    &=& \frac{1}{4} e^{-\mu (3x+2y)} \left( \frac{\lambda}{\mu+\lambda} \right)^j \cdot \left( \frac{\lambda}{\lambda-\mu} \right)^{k-j} \frac{\gamma \big( k-j, (\lambda-\mu)x \big)}{\Gamma(k-j)} \cdot \left( 2\mu x +1 - \frac{2 \mu j}{\mu+\lambda} \right)\\
    &+& \frac{1}{4} e^{-\mu (x+2y)} \cdot \left( \frac{\lambda}{\mu+\lambda} \right)^k \cdot \left( 1- \frac{\gamma \big( k-j, (\lambda+\mu)x \big)}{\Gamma(k-j)} \right) \cdot \left( 1+ \frac{2\mu j}{\mu + \lambda} \right)\\
    &+&\frac{1}{4} e^{-\mu (x+2y)} \cdot \left( \frac{\lambda}{\mu+\lambda} \right)^k \cdot \frac{k-j}{\lambda+\mu} \cdot \left( 1- \frac{\gamma \big( k-j+1, (\lambda+\mu)x \big)}{\Gamma(k-j)} \right) \\
    &-& \frac{1}{2} \mu e^{-\mu(3x+2y)} \left( \frac{\lambda}{\mu+\lambda} \right)^j \cdot \frac{j}{\lambda+\mu} \cdot \left( \frac{\lambda}{\lambda-\mu} \right)^{k-j} \cdot \frac{\gamma \big( k-j, (\lambda-\mu)x \big)}{\Gamma(k-j)}\\
    &-&\frac{1}{2} \mu e^{-\mu(x+2y)} \left( \frac{\lambda}{\mu+\lambda} \right)^k \cdot \frac{j}{\lambda+\mu} \cdot \left(1- \frac{\gamma \big( k-j, (\lambda+\mu)x \big)}{\Gamma(k-j)} \right)\\
    &-& \frac{1}{2} \mu e^{-\mu(3x+2y)} \left( \frac{\lambda}{\mu+\lambda} \right)^j \cdot \frac{k-j}{\lambda-\mu} \cdot \left( \frac{\lambda}{\lambda-\mu} \right)^{k-j} \cdot \frac{\gamma \big( k-j+1, (\lambda-\mu)x \big)}{\Gamma(k-j+1)}\\
    &-&\frac{1}{2} \mu e^{-\mu(x+2y)} \left( \frac{\lambda}{\mu+\lambda} \right)^k \cdot \frac{k-j}{\lambda+\mu} \cdot \left(1- \frac{\gamma \big( k-j+1, (\lambda+\mu)x \big)}{\Gamma(k-j+1)} \right).\\
\end{eqnarray*}

\subsubsection{$\lambda = \mu$}
Following the split of equation into following three parts in \ref{seca41},:
\begin{enumerate}
    \item \begin{eqnarray*}
        \lefteqn{ A= e^{-\mu (x+2y)} \int_{\mathbb{R}_+^3} \{z_3 > x + z_1 \} \cdot \{ z_3> z_1+z_2\}  }\\ 
        &\times& \mu^2 \cdot e^{-\mu(2z_3-z_2-z_1)} \cdot M \cdot e^{-\lambda (z_1+z_2)} \cdot z_3 dz_1 dz_2 dz_3\\
        &=& e^{-\mu (x+2y)} \int_{\mathbb{R}_+^3} \{z_3 > x + z_1 \} \cdot \{ x>z_2\} \\
        &\times& \mu^2 \cdot e^{-\mu(2z_3-z_2-z_1)} \cdot M \cdot e^{-\lambda (z_1+z_2)} \cdot z_3 dz_1 dz_2 dz_3\\
        &+& e^{-\mu (x+2y)} \int_{\mathbb{R}_+^3} \{z_3 > z_2 + z_1 \} \cdot \{ x<z_2\} \\
        &\times& \mu^2 \cdot e^{-\mu(2z_3-z_2-z_1)} \cdot M \cdot e^{-\lambda (z_1+z_2)} \cdot z_3 dz_1 dz_2 dz_3\\
        &=& \frac{1}{4} e^{-\mu (3x+2y)} \cdot (2\mu x+1) \int_{\mathbb{R}_+^2} \{ x > z_2\} \cdot e^{-(\mu+\lambda) z_1} \cdot M dz_1 dz_2 \\
        &+& \frac{1}{4} e^{-\mu (3x+2y)} \cdot (2\mu) \int_{\mathbb{R}_+^2} z_1 \cdot \{ x > z_2\} \cdot e^{-(\mu+\lambda) z_1}\cdot M dz_1 dz_2 \\
        &+& \frac{1}{4} e^{-\mu (x+2y)}  \int_{\mathbb{R}_+^2} \{ x < z_2\} \cdot e^{-(\mu+\lambda) z_2} \cdot M dz_1 dz_2 \\
        &+& \frac{1}{4} e^{-\mu (x+2y)} \cdot (2\mu) \int_{\mathbb{R}_+^2} (z_1+z_2) \cdot \{ x < z_2\} \cdot e^{-(\mu+\lambda) z_1} \cdot e^{-(\mu+\lambda) z_2} \cdot M dz_1 dz_2 \\
        &=&  \frac{1}{4} e^{-\mu (3x+2y)} \left( \frac{\lambda}{\mu+\lambda} \right)^j \cdot\frac{\lambda^{k-j}}{\Gamma(k-j)} \cdot \left( -\frac{x^{k-j}}{k-j} \right) \cdot \left( 2\mu x +1 - \frac{2 \mu j}{\mu+\lambda} \right)\\
        &+& \frac{1}{4} e^{-\mu (x+2y)} \cdot \left( \frac{\lambda}{\mu+\lambda} \right)^k \cdot \left( 1- \frac{\gamma \big( k-j, (\lambda+\mu)x \big)}{\Gamma(k-j)} \right) \cdot \left( 1+ \frac{2\mu j}{\mu + \lambda} \right)\\
        &+&\frac{1}{4} e^{-\mu (x+2y)} \cdot \left( \frac{\lambda}{\mu+\lambda} \right)^k \cdot \frac{k-j}{\lambda+\mu} \cdot \left( 1- \frac{\gamma \big( k-j+1, (\lambda+\mu)x \big)}{\Gamma(k-j)} \right). \\
    \end{eqnarray*}
        
    \item \begin{eqnarray*}
        \lefteqn{ -B= e^{-\mu (x+2y)} \int_{\mathbb{R}_+^3} \{z_3 > x + z_1 \} \cdot \{ z_3> z_1+z_2\}  }\\ 
        &\times& \mu^2 \cdot e^{-\mu(2z_3-z_2-z_1)} \cdot M \cdot e^{-\lambda (z_1+z_2)} \cdot z_1 dz_1 dz_2 dz_3 \\
        &=& e^{-\mu (x+2y)} \int_{\mathbb{R}_+^3} \{z_3 > x + z_1 \} \cdot \{ x>z_2\} \\
        &\times&\mu^2 \cdot e^{-\mu(2z_3-z_2-z_1)} \cdot M \cdot e^{-\lambda (z_1+z_2)} \cdot z_1 dz_1 dz_2 dz_3\\
        &+& e^{-\mu (x+2y)} \int_{\mathbb{R}_+^3} \{z_3 > z_2 + z_1 \} \cdot \{ x<z_2\} \\
        &\times&\mu^2 \cdot e^{-\mu(2z_3-z_2-z_1)} \cdot M \cdot e^{-\lambda (z_1+z_2)} \cdot z_1 dz_1 dz_2 dz_3\\
        &=& \frac{1}{2} \mu e^{-\mu(3x+2y)}  \int_{\mathbb{R}_+^2} \{x > z_2\} \cdot e^{-(\lambda+\mu)z_1}  \cdot M \cdot z_1 dz_1 dz_2\\
        &+& \frac{1}{2} \mu e^{-\mu(x+2y)}  \int_{\mathbb{R}_+^2} \{x < z_2\} \cdot e^{-(\lambda+\mu)z_1} \cdot e^{-(\lambda+\mu)z_2} \cdot M \cdot z_1 dz_1 dz_2\\
        &=& \frac{1}{2} \mu e^{-\mu(3x+2y)} \left( \frac{\lambda}{\mu+\lambda} \right)^j \cdot \frac{j}{\lambda+\mu} \cdot \frac{\lambda^{k-j}}{\Gamma(k-j)} \cdot \left(\frac{x^{k-j}}{k-j} \right) \\
        &+&\frac{1}{2} \mu e^{-\mu(x+2y)} \left( \frac{\lambda}{\mu+\lambda} \right)^k \cdot \frac{j}{\lambda+\mu} \cdot \left(1- \frac{\gamma \big( k-j, (\lambda+\mu)x \big)}{\Gamma(k-j)} \right).\\
    \end{eqnarray*}

    \item \begin{eqnarray*}
        \lefteqn{ -C= e^{-\mu (x+2y)} \int_{\mathbb{R}_+^3} \{z_3 > x + z_1 \} \cdot \{ z_3> z_1+z_2\} }\\ 
        &\times& \mu^2 \cdot e^{-\mu(2z_3-z_2-z_1)} \cdot M \cdot e^{-\lambda (z_1+z_2)} \cdot z_2 dz_1 dz_2 dz_3\\
        &=& e^{-\mu (x+2y)} \int_{\mathbb{R}_+^3} \{z_3 > x + z_1 \} \cdot \{ x>z_2\} \\
        &\times&\mu^2 \cdot e^{-\mu(2z_3-z_2-z_1)} \cdot M \cdot e^{-\lambda (z_1+z_2)} \cdot z_2 dz_1 dz_2 dz_3\\
        &+& e^{-\mu (x+2y)} \int_{\mathbb{R}_+^3} \{z_3 > z_2 + z_1 \} \cdot \{ x<z_2\} \\
        &\times& \mu^2 \cdot e^{-\mu(2z_3-z_2-z_1)} \cdot M \cdot e^{-\lambda (z_1+z_2)} \cdot z_2 dz_1 dz_2 dz_3\\
        &=& \frac{1}{2} \mu e^{-\mu(3x+2y)} \left( \frac{\lambda}{\mu+\lambda} \right)^j \cdot \frac{\lambda^{k-j}}{\Gamma(k-j)} \cdot \frac{x^{k-j+1}}{k-j+1}\\
        &+&\frac{1}{2} \mu e^{-\mu(x+2y)} \left( \frac{\lambda}{\mu+\lambda} \right)^j \cdot \frac{k-j}{\lambda+\mu} \cdot \left(1- \frac{\gamma \big( k-j+1, (\lambda+\mu)x \big)}{\Gamma(k-j+1)} \right).\\
    \end{eqnarray*}
\end{enumerate}
Putting all three parts together we have:
\begin{eqnarray*}
    \lefteqn{ \mathbb{P} \left( S^{(1)}_n - \mathcal{A}_{j} > x , S^{(1)}_{n+j} > x, S^{(2)}_{n+k} > y, S^{(1)}_{n} + S^{(2)}_{n} - \mathcal{A}_{j} - \Tilde{\mathcal{A}}_{k-j} - S^{(1)}_{n+k} > y , \mathcal{E}^C_{(2)}  \right) } \\ 
    &=& A-B-C \\
    &=& \frac{1}{4} e^{-\mu (3x+2y)} \left( \frac{\lambda}{\mu+\lambda} \right)^j \cdot\frac{\lambda^{k-j}}{\Gamma(k-j)} \cdot \left( -\frac{x^{k-j}}{k-j} \right) \cdot \left( 2\mu x +1 - \frac{2 \mu j}{\mu+\lambda} \right)\\
    &+& \frac{1}{4} e^{-\mu (x+2y)} \cdot \left( \frac{\lambda}{\mu+\lambda} \right)^k \cdot \left( 1- \frac{\gamma \big( k-j, (\lambda+\mu)x \big)}{\Gamma(k-j)} \right) \cdot \left( 1+ \frac{2\mu j}{\mu + \lambda} \right)\\
    &+&\frac{1}{4} e^{-\mu (x+2y)} \cdot \left( \frac{\lambda}{\mu+\lambda} \right)^k \cdot \frac{k-j}{\lambda+\mu} \cdot \left( 1- \frac{\gamma \big( k-j+1, (\lambda+\mu)x \big)}{\Gamma(k-j)} \right) \\
    &-& \frac{1}{2} \mu e^{-\mu(3x+2y)} \left( \frac{\lambda}{\mu+\lambda} \right)^j \cdot \frac{j}{\lambda+\mu} \cdot \frac{\lambda^{k-j}}{\Gamma(k-j)} \cdot \left(\frac{x^{k-j}}{k-j} \right) \\
    &-&\frac{1}{2} \mu e^{-\mu(x+2y)} \left( \frac{\lambda}{\mu+\lambda} \right)^k \cdot \frac{j}{\lambda+\mu} \cdot \left(1- \frac{\gamma \big( k-j, (\lambda+\mu)x \big)}{\Gamma(k-j)} \right)\\
    &-& \frac{1}{2} \mu e^{-\mu(3x+2y)} \left( \frac{\lambda}{\mu+\lambda} \right)^j \cdot \frac{\lambda^{k-j}}{\Gamma(k-j)} \cdot \frac{x^{k-j+1}}{k-j+1}\\
    &-&\frac{1}{2} \mu e^{-\mu(x+2y)} \left( \frac{\lambda}{\mu+\lambda} \right)^j \cdot \frac{k-j}{\lambda+\mu} \cdot \left(1- \frac{\gamma \big( k-j+1, (\lambda+\mu)x \big)}{\Gamma(k-j+1)} \right).\\
\end{eqnarray*}

\subsubsection{$\lambda<\mu$}
Applying similar technique seen in \ref{A33} We have:
\begin{eqnarray*}
    \lefteqn{ \mathbb{P} \left( S^{(1)}_n - \mathcal{A}_{j} > x , S^{(1)}_{n+j} > x, S^{(2)}_{n+k} > y, S^{(1)}_{n} + S^{(2)}_{n} - \mathcal{A}_{j} - \Tilde{\mathcal{A}}_{k-j} - S^{(1)}_{n+k} > y , \mathcal{E}^C_{(2)}  \right) } \\ 
    &=& A-B-C \\
    &=& \frac{1}{4} e^{-\mu (3x+2y)} \left( \frac{\lambda}{\mu+\lambda} \right)^j \cdot \left( \frac{\lambda}{\lambda-\mu} \right)^{k-j} \frac{-\gamma \big( k-j, (-\lambda+\mu)x \big)}{\Gamma(k-j)} \cdot \left( 2\mu x +1 - \frac{2 \mu j}{\mu+\lambda} \right)\\
    &+& \frac{1}{4} e^{-\mu (x+2y)} \cdot \left( \frac{\lambda}{\mu+\lambda} \right)^k \cdot \left( 1- \frac{\gamma \big( k-j, (\lambda+\mu)x \big)}{\Gamma(k-j)} \right) \cdot \left( 1+ \frac{2\mu j}{\mu + \lambda} \right)\\
    &+&\frac{1}{4} e^{-\mu (x+2y)} \cdot \left( \frac{\lambda}{\mu+\lambda} \right)^k \cdot \frac{k-j}{\lambda+\mu} \cdot \left( 1- \frac{\gamma \big( k-j+1, (\lambda+\mu)x \big)}{\Gamma(k-j)} \right) \\
    &-& \frac{1}{2} \mu e^{-\mu(3x+2y)} \left( \frac{\lambda}{\mu+\lambda} \right)^j \cdot \frac{j}{\lambda+\mu} \cdot \left( \frac{\lambda}{\lambda-\mu} \right)^{k-j} \cdot \frac{-\gamma \big( k-j, (-\lambda+\mu)x \big)}{\Gamma(k-j)}\\
    &-&\frac{1}{2} \mu e^{-\mu(x+2y)} \left( \frac{\lambda}{\mu+\lambda} \right)^k \cdot \frac{j}{\lambda+\mu} \cdot \left(1- \frac{\gamma \big( k-j, (\lambda+\mu)x \big)}{\Gamma(k-j)} \right)\\
    &-& \frac{1}{2} \mu e^{-\mu(3x+2y)} \left( \frac{\lambda}{\mu+\lambda} \right)^j \frac{\lambda^{(k-j)}}{\Gamma(k-j) } \\
&\times& \left(  \sum_{i=0}^{k-j-1} (-1)^{i} \cdot \frac{e^{(\mu-\lambda)x}(k-j-1)!}{(k-j-1-i)!\cdot \mu-\lambda)^{i+1}} \cdot x^{k-j-1-i} - (-1)^{k-j-1}\frac{(k-j-1)!}{\mu-\lambda)^{k-j-1}} \right)\\
    &-&\frac{1}{2} \mu e^{-\mu(x+2y)} \left( \frac{\lambda}{\mu+\lambda} \right)^k \cdot \frac{k-j}{\lambda+\mu} \cdot \left(1- \frac{\gamma \big( k-j+1, (\lambda-\mu)x \big)}{\Gamma(k-j+1)} \right).
\end{eqnarray*}

\section{Proof for CASE IV in \ref{ovlp_diff_customers_prop}}
\subsection{Condition 1:}
Let $\mathcal{A}_{j}, \mathcal{A}_{m}, \mathcal{A}_{k}$ be the inter-arrival times between customer $n$ and $n+j$, \ $n+j$ and $m$, $m$ and $m+k$ respectively. In addition, we define $M$ to be $$ M = \frac{\lambda^{(j)+(k)} z_1^{j-1} z_2^{k-1}}{\Gamma(j)  \Gamma(k)}.$$This implies that:
\begin{eqnarray*}
    \lefteqn{ \mathbb{P} \left( S^{(1)}_n - \mathcal{A}_{j} > x , S^{(1)}_{n+j} > x, S^{(2)}_{m+k} > y, S^{(1)}_{m} + S^{(2)}_{m} - \mathcal{A}_{k} - S^{(1)}_{m+k} > y , \mathcal{E}_{(2)}  \right) } \\  
    &=& \int_{\mathbb{R}_+^8} \{z_3 > x + z_1\} \cdot \{ z_4 > x \} \cdot \{ z_8 > y \} \cdot \{ z_5 + z_6 - z_2 - z_7 > y \} \cdot \{ z_7 + z_2 > z_5 \}\\
    &\times& \mu^6 e^{-\mu(z_3+z_4+z_5+z_6+z_7+z_8)} \cdot M \cdot e^{-\lambda (z_1+z_2)} dz_1 dz_2 dz_3 dz_4 dz_5 dz_6 dz_7 dz_8 \\
    &=& e^{-\mu(2x+y)} \int_{\mathbb{R}_+^5} \{ z_5 + z_6 - z_2 - z_7 > y \} \cdot \{ z_7 + z_2 > z_5 \} \cdot \mu^3 e^{-\mu(z_1+z_5+z_6+z_7)}  \\
    &\times& M \cdot e^{-\lambda (z_1+z_2)} dz_1 dz_2 dz_5 dz_6 dz_7 \\
    &=& e^{-\mu(2x+2y)} \int_{\mathbb{R}_+^4}  \{ z_7 + z_2 > z_5 \} \cdot \mu^2 e^{-\mu(z_1+z_2+2z_7)} \cdot M \cdot e^{-\lambda (z_1+z_2)} dz_1 dz_2 dz_5 dz_7 \\
    &=& \frac{1}{2} e^{-\mu(2x+2y)} \int_{\mathbb{R}_+^3} \{z_5>z_2\} \cdot \mu e^{-\mu(z_1-z_2+2z_5)} \cdot M \cdot e^{-\lambda (z_1+z_2)} dz_1 dz_2 dz_5 \\
    &+& \frac{1}{2} e^{-\mu(2x+2y)} \int_{\mathbb{R}_+^3}  \{ z_2 > z_5 \} \cdot \mu e^{-\mu(z_1+z_2)} \cdot M \cdot e^{-\lambda (z_1+z_2)} dz_1 dz_2 dz_5 \\
    &=& \frac{1}{4} e^{-\mu(2x+2y)} \int_{\mathbb{R}_+^3} e^{-(\mu+\lambda)z_1} e^{-(\mu+\lambda)z_2}  \cdot M dz_1 dz_2 \\
    &+& \frac{1}{2} \mu e^{-\mu(2x+2y)} \int_{\mathbb{R}_+^3}  z_2  \cdot e^{-(\mu+\lambda)z_1} e^{-(\mu+\lambda)z_2} \cdot M  dz_1 dz_2 \\
    &=& \frac{1}{4} e^{-\mu(2x+2y)} \left( \frac{\lambda}{\lambda+\mu} \right)^{j+k}+ \frac{1}{2} \mu e^{-\mu(2x+2y)} \left( \frac{\lambda}{\lambda+\mu} \right)^{j+k}  \frac{k}{\lambda+\mu}\\
    &=&\frac{1}{4} e^{-\mu(2x+2y)} \left( \frac{\lambda}{\lambda+\mu} \right)^{j+k} \cdot \left( 1 + \frac{2\mu k}{\lambda+\mu} \right).
\end{eqnarray*}

\subsection{Condition 2:}
Let $\mathcal{A}_{j}, \mathcal{A}_{m}, \mathcal{A}_{k}$ be the inter-arrival times between customer $n$ and $n+j$, \ $n+j$ and $m$, $m$ and $m+k$ respectively. In addition, we define $M$ to be  $$ M = \frac{\lambda^{(j)+(k)} z_1^{j-1} z_2^{k-1}}{\Gamma(j)  \Gamma(k)}.$$This implies that:
\begin{eqnarray*}
   \lefteqn{ \mathbb{P} \left( S^{(1)}_n - \mathcal{A}_{j} > x , S^{(1)}_{n+j} > x, S^{(2)}_{m} > y, S^{(1)}_{m+k} + S^{(2)}_{m+k} + \mathcal{A}_{k} - S^{(1)}_{m} > y , \mathcal{E}_{(2)}^C  \right) }\\  
    &=& \int_{\mathbb{R}_+^8} \{z_3 > x + z_1\} \cdot \{ z_4 > x \} \cdot \{ z_6 > y \} \cdot \{ -z_5 + z_8 + z_2 + z_7 > y \} \cdot \{ z_7 + z_2 < z_5 \}\\
    &\times& \mu^6 e^{-\mu(z_3+z_4+z_5+z_6+z_7+z_8)} \cdot M \cdot e^{-\lambda (z_1+z_2)} dz_1 dz_2 dz_3 dz_4 dz_5 dz_6 dz_7 dz_8 \\
    &=& e^{-\mu(2x+y)} \int_{\mathbb{R}_+^5} \{ -z_5 + z_8 + z_2 + z_7 > y \} \cdot \{ z_7 + z_2 < z_5 \}\\
    &\times& \mu^3 e^{-\mu(z_1+z_5+z_7+z_8)} \cdot M \cdot e^{-\lambda (z_1+z_2)} dz_1 dz_2 dz_5 dz_7 dz_8 \\
    &=& e^{-\mu(2x+2y)} \int_{\mathbb{R}_+^4} \{ z_7 + z_2 < z_5 \} \cdot \mu^2 e^{-\mu(z_1-z_2+2z_5)} \cdot M \cdot e^{-\lambda (z_1+z_2)} dz_1 dz_2 dz_5 dz_7 \\
    &=& \frac{1}{2} e^{-\mu(2x+2y)} \int_{\mathbb{R}_+^3}  \mu e^{-\mu(z_1+z_2+2z_7)} \cdot M \cdot e^{-\lambda (z_1+z_2)} dz_1 dz_2 dz_7 \\
    &=& \frac{1}{4} e^{-\mu(2x+2y)} \int_{\mathbb{R}_+^2} e^{-\mu(z_1+z_2)} \cdot M \cdot e^{-\lambda (z_1+z_2)} dz_1 dz_2 \\
    &=& \frac{1}{4} e^{-\mu(2x+2y)}  \left( \frac{\lambda}{\mu + \lambda} \right)^{k+j}.
\end{eqnarray*}

\section{Proof for CASE V in \ref{ovlp_diff_customers_prop}}
\subsection{Condition 1:}
Let $\mathcal{A}_{m}, \mathcal{A}_{k}, \mathcal{A}_{j}$ be the inter-arrival times between customer $n$ and $m$, \ $m$ and $m+k$, $m+k$ and $n+j$ respectively. In addition, we define $M$ to be $$ M = \frac{\lambda^{(m-n) + (k) + (n+j-m-k)} z_1^{m-n-1} z_2^{k-1}z_3^{n+j-m-k}}  {\Gamma(j)  \Gamma(k) \Gamma(n+j-m-k)}.$$This implies that:
\begin{eqnarray*}
    \lefteqn{ \mathbb{P} \left( S^{(1)}_n - \mathcal{A}_{j} - \mathcal{A}_{k} - \mathcal{A}_{m} > x , S^{(1)}_{n+j} > x, S^{(2)}_{m+k} > y, S^{(1)}_{m} + S^{(2)}_{m} - \mathcal{A}_{k} - S^{(1)}_{m+k} > y , \mathcal{E}_{(2)}  \right) } \\  
    &=& \int_{\mathbb{R}_+^9} \{z_4 > x + z_1 + z_2 + z_3\} \cdot \{ z_7 > x \} \cdot  \{ z_5 + z_8 - z_2 - z_6 > y \} \cdot \{ z_6 + z_2 > z_5 \}\\
    &\times& \{ z_9 > y \} \cdot \mu^6 e^{-\mu(z_4+z_5+z_6+z_7+z_8+z_9)} \cdot M \cdot e^{-\lambda (z_1+z_2+z_3)} dz_1 dz_2 dz_3 dz_4 dz_5 dz_6 dz_7 dz_8 dz_9\\
    &=& e^{-\mu(2x+y)}\int_{\mathbb{R}_+^6} \{ z_5 + z_8 - z_2 - z_6 > y \} \cdot \{ z_6 + z_2 > z_5 \}\\
    &\times& \mu^3 e^{-\mu(z_1+z_2+z_3+z_5+z_6+z_8)} \cdot M \cdot e^{-\lambda (z_1+z_2+z_3)} dz_1 dz_2 dz_3 dz_5 dz_6 dz_8\\
    &=& e^{-\mu(2x+2y)}\int_{\mathbb{R}_+^5} \{ z_6 + z_2 > z_5 \} \cdot \mu^2 e^{-\mu(z_1+2z_2+z_3+2z_6)} \cdot M \cdot e^{-\lambda (z_1+z_2+z_3)} dz_1 dz_2 dz_3 dz_5 dz_6\\
    &=& \frac{1}{2} e^{-\mu(2x+2y)}\int_{\mathbb{R}_+^4} \{z_5 > z_2\} \mu e^{-\mu(z_1+z_3+2z_5)} \cdot M \cdot e^{-\lambda (z_1+z_2+z_3)} dz_1 dz_2 dz_3 dz_5\\
    &+& \frac{1}{2} e^{-\mu(2x+2y)}\int_{\mathbb{R}_+^4} \{ z_5 < z_2 \} \cdot \mu e^{-\mu(z_1+2z_2+z_3)} \cdot M \cdot e^{-\lambda (z_1+z_2+z_3)} dz_1 dz_2 dz_3 dz_5\\
    &=& \frac{1}{4} e^{-\mu(2x+2y)}\int_{\mathbb{R}_+^3} e^{-\mu(z_1+2z_2+z_3)} \cdot M \cdot e^{-\lambda (z_1+z_2+z_3)} dz_1 dz_2 dz_3\\
    &+& \frac{1}{2} \mu e^{-\mu(2x+2y)}\int_{\mathbb{R}_+^3} z_2 \cdot e^{-\mu(z_1+2z_2+z_3)} \cdot M \cdot e^{-\lambda (z_1+z_2+z_3)} dz_1 dz_2 dz_3\\
    &=& \frac{1}{4} e^{-\mu(2x+2y)} \left( \frac{\lambda}{2\mu + \lambda} \right)^{k} \left( \frac{\lambda}{\mu + \lambda} \right)^{j-k} \cdot \left( 1+ \frac{2\mu k}{2\mu + \lambda} \right).
\end{eqnarray*}

\subsection{Condition 2:}
Let $\mathcal{A}_{m}, \mathcal{A}_{k}, \mathcal{A}_{j}$ be the inter-arrival times between customer $n$ and $m$, \ $m$ and $m+k$, $m+k$ and $n+j$ respectively. In addition, we define $M$ to be $$ M = \frac{\lambda^{(m-n) + (k) + (n+j-m-k)} z_1^{m-n-1} z_2^{k-1}z_3^{n+j-m-k}}  {\Gamma(j)  \Gamma(k) \Gamma(n+j-m-k)}.$$This implies that:
\begin{eqnarray*}
   \lefteqn{ \mathbb{P} \left( S^{(1)}_n - \mathcal{A}_{j} - \mathcal{A}_{k} - \mathcal{A}_{m} > x , S^{(1)}_{n+j} > x, S^{(2)}_{m} > y, S^{(1)}_{m+k} + S^{(2)}_{m+k} + \mathcal{A}_{k} - S^{(1)}_{m} > y , \mathcal{E}_{(2)}^C  \right) } \\  
    &=& \int_{\mathbb{R}_+^9}\{z_4 > x + z_1 + z_2 + z_3\} \cdot \{ z_7 > x \}  \cdot \{ -z_5 + z_9 + z_2 + z_6 > y \} \cdot \{ z_6 + z_2 < z_5 \}\\
    &\times&  \{ z_8 > y \} \cdot \mu^6 e^{-\mu(z_4+z_5+z_6+z_7+z_8+z_9)} \cdot M \cdot e^{-\lambda (z_1+z_2+z_3)} dz_1 dz_2 dz_3 dz_4 dz_5 dz_6 dz_7 dz_8 dz_9\\
    &=& e^{-\mu(2x+y)}\int_{\mathbb{R}_+^6}  \{ -z_5 + z_9 + z_2 + z_6 > y \} \cdot \{ z_6 + z_2 < z_5 \}\\
    &\times& \mu^3 e^{-\mu(z_1+z_2+z_3+z_5+z_6+z_9)} \cdot M \cdot e^{-\lambda (z_1+z_2+z_3)} dz_1 dz_2 dz_3 dz_5 dz_6 dz_9\\
    &=& \frac{1}{2} e^{-\mu(2x+2y)}\int_{\mathbb{R}_+^4}\mu e^{-\mu(z_1+2z_2+z_3+2z_6)} \cdot M \cdot e^{-\lambda (z_1+z_2+z_3)} dz_1 dz_2 dz_3  dz_6 \\
    &=& \frac{1}{4}  e^{-\mu(2x+2y)}\int_{\mathbb{R}_+^3} e^{-\mu(z_1+2z_2+z_3)} \cdot M \cdot e^{-\lambda (z_1+z_2+z_3)} dz_1 dz_2 dz_3 \\
    &=& \frac{1}{4} e^{-\mu(2x+2y)} \left( \frac{\lambda}{2\mu + \lambda} \right)^{k} \left( \frac{\lambda}{\mu + \lambda} \right)^{j-k}.
\end{eqnarray*}

\section{Proof for CASE VI in \ref{ovlp_diff_customers_prop}}
\subsection{Condition 1:}
Let $\mathcal{A}_{k}, \mathcal{A}_{j-k}$ be the inter-arrival times between customer $n$ and $n+k$, $n+k$ and $n+j$ respectively. In addition, we define $M$ to be $$ M = \frac{\lambda^{(k) + (j-k)} z_1^{k-1} z_2^{j-k-1}}  {\Gamma(k) \Gamma(j-k)}.$$This implies that:
\begin{eqnarray*}
    \lefteqn{ \mathbb{P} \left( S^{(1)}_n - \mathcal{A}_{k} - \mathcal{A}_{j-k} > x , S^{(1)}_{n+j} > x, S^{(2)}_{n+k} > y, S^{(1)}_{n} + S^{(2)}_{n} - \mathcal{A}_{k} - S^{(1)}_{n+k} > y , \mathcal{E}_{(2)}  \right) } \\  
    &=& \int_{\mathbb{R}_+^7} \{z_3 > x + z_1 + z_2\} \cdot \{ z_5 > x \} \cdot \{ z_7 > y \} \cdot \{ z_3 + z_6 - z_1 - z_4 > y \} \cdot \{ z_1 + z_4 > z_3 \}\\
    &\times& \mu^5 e^{-\mu(z_3+z_4+z_5+z_6+z_7)} \cdot M \cdot e^{-\lambda (z_1+z_2)} dz_1 dz_2 dz_3 dz_4 dz_5 dz_6 dz_7\\
    &=& e^{-\mu(x+2y)} \int_{\mathbb{R}_+^4} \{z_3 > x + z_1 + z_2\} \cdot \{ z_1 + z_4 > z_3 \} \\
    &\times&\mu^2 e^{-\mu(z_1+2z_4)} \cdot M \cdot e^{-\lambda (z_1+z_2)} dz_1 dz_2 dz_3 dz_4 \\
    &=& \frac{1}{2} e^{-\mu(x+2y)} \int_{\mathbb{R}_+^3} \{z_3 > x + z_1 + z_2\} \cdot \mu e^{-\mu(-z_1+2z_3)} \cdot M \cdot e^{-\lambda (z_1+z_2)} dz_1 dz_2 dz_3 \\
    &=& \frac{1}{4} e^{-\mu(3x+2y)} \int_{\mathbb{R}_+^2}  e^{-\mu(z_1+2z_2)} \cdot M \cdot e^{-\lambda (z_1+z_2)} dz_1 dz_2 \\
    &=& \frac{1}{4} e^{-\mu(3x+2y)}  \left( \frac{\lambda}{2\mu + \lambda} \right)^{j-k} \left( \frac{\lambda}{\mu + \lambda} \right)^{k}.
\end{eqnarray*}

\subsection{Condition 2:}
Let $\mathcal{A}_{k}, \mathcal{A}_{j-k}$ be the inter-arrival times between customer $n$ and $n+k$, $n+k$ and $n+j$ respectively. In addition, we define $M$ to be $$ M = \frac{\lambda^{(k) + (j-k)} z_1^{k-1} z_2^{j-k-1}}  {\Gamma(k) \Gamma(j-k)}.$$This implies that:
\begin{eqnarray*}
   \lefteqn{ \mathbb{P} \left( S^{(1)}_n - \mathcal{A}_{k} - \mathcal{A}_{j-k} > x , S^{(1)}_{n+j} > x, S^{(2)}_{n} > y, S^{(1)}_{n+k} + S^{(2)}_{n+k} + \mathcal{A}_{k} - S^{(1)}_{n} > y , \mathcal{E}_{(2)}^C  \right) } \\ 
    &=& \int_{\mathbb{R}_+^7} \{z_3 > x + z_1 + z_2\} \cdot \{ z_5 > x \} \cdot \{ z_6 > y \} \cdot \{ -z_3 + z_7 + z_1 + z_4 > y \}  \\
    &\times& \{ z_1 + z_4 < z_3 \} \cdot \mu^5 e^{-\mu(z_3+z_4+z_5+z_6+z_7)} \cdot M \cdot e^{-\lambda (z_1+z_2)} dz_1 dz_2 dz_3 dz_4 dz_5 dz_6 dz_7\\
    &=& e^{-\mu(x+2y)} \int_{\mathbb{R}_+^4}  \{z_3 > x + z_1 + z_2\} \cdot\{ z_1 + z_4 < z_3 \} \\
    &\times&\mu^2 e^{-\mu(-z_1+2z_3)} \cdot M \cdot e^{-\lambda (z_1+z_2)} dz_1 dz_2 dz_3 dz_4 \\
    &=& e^{-\mu(x+2y)} \int_{\mathbb{R}_+^4}  \{z_3 > x + z_1 + z_2\} \cdot z_3 \cdot \mu^2 e^{-\mu(-z_1+2z_3)} \cdot M \cdot e^{-\lambda (z_1+z_2)} dz_1 dz_2 dz_3 \\
    &-& e^{-\mu(x+2y)} \int_{\mathbb{R}_+^4}  \{z_3 > x + z_1 + z_2\} \cdot z_1 \cdot \mu^2 e^{-\mu(-z_1+2z_3)} \cdot M \cdot e^{-\lambda (z_1+z_2)} dz_1 dz_2 dz_3 \\
    &=& \frac{1}{4} e^{-\mu(3x+2y)}\int_{\mathbb{R}_+^3} e^{-\mu(z_1+2z_2)} \cdot (2\mu x +1 +2\mu z_1+ 2\mu z_2) \cdot M \cdot e^{-\lambda (z_1+z_2)} dz_1 dz_2 \\
    &-& \frac{1}{2} \mu e^{-\mu(3x+2y)} \int_{\mathbb{R}_+^3} z_1 \cdot \mu e^{-\mu(z_1+2z_2)} \cdot M \cdot e^{-\lambda (z_1+z_2)} dz_1 dz_2 \\
    &=& \frac{1}{4} e^{-\mu(3x+2y)}  \left( \frac{\lambda}{2\mu + \lambda} \right)^{j-k} \left( \frac{\lambda}{\mu + \lambda} \right)^{k} \cdot \left( 2\mu x +1 + 2\mu \cdot \frac{j-k}{2\mu+\lambda} \right).
\end{eqnarray*}

\section{Proof for CASE VII in \ref{ovlp_diff_customers_prop}}
\subsection{Condition 1:}
Let $\mathcal{A}_{n}, \mathcal{A}_{k}, \mathcal{A}_{j}$ be the inter-arrival times between customer $m$, $n$, $n$ and $m+k$, $m+k$ and $n+j$ respectively. In addition, we define $M$ to be $$ M = \frac{\lambda^{(n-m) + (m-n+k) + (n+j-m-k)} z_1^{n-m-1} z_2^{m-n+k-1} z_3^{n+j-m-k-1}}  {\Gamma(n-m) \Gamma(m-n+k) \Gamma(n+j-m-k)}.$$This implies that:
\begin{eqnarray*}
    \lefteqn{ \mathbb{P} \left( S^{(1)}_n - \mathcal{A}_{k} - \mathcal{A}_{j} > x , S^{(1)}_{n+j} > x, S^{(2)}_{m+k} > y, S^{(1)}_{m} + S^{(2)}_{m} - \mathcal{A}_{k} - \mathcal{A}_{n} - S^{(1)}_{m+k} > y , \mathcal{E}_{(2)}  \right) } \\  
    &=& \int_{\mathbb{R}_+^9} \{z_4 > x + z_3 + z_2\} \cdot \{ z_6 > x \} \cdot\{ z_9 + z_7 - z_1 - z_2 - z_5 > y \} \\
    &\times& \{ z_1 + z_2 + z_5 > z_9 \} \cdot \{ z_8 > y \} \cdot \mu^6 e^{-\mu(z_4+z_5+z_6+z_7+z_8+z_9)} \\
    &\times& M \cdot e^{-\lambda (z_1+z_2+z_3)} dz_1 dz_2 dz_3 dz_4 dz_5 dz_6 dz_7 dz_8 dz_9\\
    &=& e^{-\mu(2x+2y)} \int_{\mathbb{R}_+^5}  \{ z_1 + z_2 + z_5 > z_9 \} \\
    &\times& \mu^2 e^{-\mu(z_1+2z_2+z_3+2z_5)} \cdot M \cdot e^{-\lambda (z_1+z_2+z_3)} dz_1 dz_2 dz_3 dz_5 dz_9\\
    &=& \frac{1}{2} e^{-\mu(2x+2y)} \int_{\mathbb{R}_+^4}  \{ z_9 > z_1+z_2 \} \cdot \mu e^{-\mu(-z_1+z_3+2z_9)} \cdot M \cdot e^{-\lambda (z_1+z_2+z_3)} dz_1 dz_2 dz_3 dz_9\\
    &+& \frac{1}{2} e^{-\mu(2x+2y)} \int_{\mathbb{R}_+^4} \{ z_9 < z_1+z_2 \} \cdot \mu e^{-\mu(z_1+2z_2+z_3)} \cdot M \cdot e^{-\lambda (z_1+z_2+z_3)} dz_1 dz_2 dz_3 dz_9\\
    &=& \frac{1}{4} e^{-\mu(2x+2y)}\int_{\mathbb{R}_+^3} e^{-(\mu+\lambda)z_1} e^{-(2\mu+\lambda)z_2} e^{-(\mu+\lambda)z_3} \cdot Mdz_1 dz_2 dz_3\\
    &+& \frac{1}{2} \mu e^{-\mu(2x+2y)} \int_{\mathbb{R}_+^3} ( z_1+z_2 ) \cdot e^{-(\mu+\lambda)z_1} e^{-(2\mu+\lambda)z_2} e^{-(\mu+\lambda)z_3} \cdot Mdz_1 dz_2 dz_3\\
    &=&\frac{1}{4} e^{-\mu(2x+2y)} \left( \frac{\lambda}{\mu + \lambda} \right)^{2n+j-2m-k} \left( \frac{\lambda}{2\mu + \lambda} \right)^{m-n+k} \cdot \left( 1+ 2\mu\left( \frac{n-m}{\lambda+\mu} + \frac{m-n+k}{2\mu + \lambda} \right) \right).
\end{eqnarray*}

\subsection{Condition 2:}
Let $\mathcal{A}_{n}, \mathcal{A}_{k}, \mathcal{A}_{j}$ be the inter-arrival times between customer $m$, $n$, $n$ and $m+k$, $m+k$ and $n+j$ respectively. In addition, we define $M$ to be $$ M = \frac{\lambda^{(n-m) + (m-n+k) + (n+j-m-k)} z_1^{n-m-1} z_2^{m-n+k-1} z_3^{n+j-m-k-1}}  {\Gamma(n-m) \Gamma(m-n+k) \Gamma(n+j-m-k)}.$$This implies that:
\begin{eqnarray*}
    \lefteqn{ \mathbb{P} \left( S^{(1)}_n - \mathcal{A}_{k} - \mathcal{A}_{j} > x , S^{(1)}_{n+j} > x, S^{(2)}_{m} > y, -S^{(1)}_{m} + S^{(2)}_{m+k} + \mathcal{A}_{k} + \mathcal{A}_{n} + S^{(1)}_{m+k} > y , \mathcal{E}_{(2)}^C  \right) } \\  
    &=& \int_{\mathbb{R}_+^9} \{z_4 > x + z_3 + z_2\} \cdot \{ z_6 > x \} \cdot \{ -z_9 + z_8 + z_1 + z_2 + z_5 > y \} \\
    &\times& \{ z_1 + z_2 + z_5 < z_9 \} \cdot \{ z_7 > y \} \cdot \mu^6 e^{-\mu(z_4+z_5+z_6+z_7+z_8+z_9)} \\
    &\times& M \cdot e^{-\lambda (z_1+z_2+z_3)} dz_1 dz_2 dz_3 dz_4 dz_5 dz_6 dz_7 dz_8 dz_9\\
    &=& e^{-\mu(2x+2y)} \int_{\mathbb{R}_+^5} \{ z_1 + z_2 + z_5 < z_9 \} \\
    &\times& \mu^2 e^{-\mu(-z_1+z_3+2z_9)} \cdot M \cdot e^{-\lambda (z_1+z_2+z_3)} dz_1 dz_2 dz_3 dz_5 dz_9 \\
    &=& \frac{1}{4} e^{-\mu(2x+2y)} \int_{\mathbb{R}_+^3}  e^{-\mu(z_1+2z_2+z_3)} \cdot M \cdot e^{-\lambda (z_1+z_2+z_3)} dz_1 dz_2 dz_3 \\
    &=&\frac{1}{4} e^{-\mu(2x+2y)} \left( \frac{\lambda}{\mu + \lambda} \right)^{2n+j-2m-k} \left( \frac{\lambda}{2\mu + \lambda} \right)^{m-n+k}.
\end{eqnarray*}

\section{Proof for CASE VIII in \ref{ovlp_diff_customers_prop}}
\subsection{Condition 1:}
Let $\mathcal{A}_{k}, \mathcal{A}_{j}$ be the inter-arrival times between customer $m$, $m+k$, $n$ and $n+j$ respectively. In addition, we define $M$ to be $$ M = \frac{\lambda^{(k) + (j)} z_1^{k-1} z_2^{j-1}}  {\Gamma(k) \Gamma(j)}.$$This implies that:
\begin{eqnarray*}
    \lefteqn{ \mathbb{P} \left( S^{(1)}_n - \mathcal{A}_{j} > x , S^{(1)}_{n+j} > x, S^{(2)}_{m+k} > y, S^{(1)}_{m} + S^{(2)}_{m} - \mathcal{A}_{k} - S^{(1)}_{m+k} > y , \mathcal{E}_{(2)}  \right) } \\  
    &=& \int_{\mathbb{R}_+^8} \{z_4 > x + z_3\} \cdot \{ z_5 > x \} \cdot \{ z_9 > y \} \cdot \{ z_6 + z_7 - z_1 - z_8 > y \} \cdot \{ z_1 + z_8 > z_6 \}\\
    &\times& \mu^6 e^{-\mu(z_4+z_5+z_6+z_7+z_8+z_9)} \cdot M \cdot e^{-\lambda (z_1+z_3)} dz_1 dz_3 dz_4 dz_5 dz_6 dz_7 dz_8 dz_9\\
    &=& e^{-\mu(2x+2y)} \int_{\mathbb{R}_+^4} \{ z_1 + z_8 > z_6 \} \cdot \mu^2 e^{-\mu(z_1+z_3+2z_8)} \cdot M \cdot e^{-\lambda (z_1+z_3)} dz_1 dz_3 dz_6 dz_8\\ 
    &=& \frac{1}{2} e^{-\mu(2x+2y)} \int_{\mathbb{R}_+^3} \{z_6 > z_1 \} \cdot \mu e^{-\mu(-z_1+z_3+2z_6)} \cdot M \cdot e^{-\lambda (z_1+z_3)} dz_1 dz_3 dz_6\\ 
    &+& \frac{1}{2} e^{-\mu(2x+2y)} \int_{\mathbb{R}_+^3} \{z_6 < z_1 \} \cdot \mu e^{-\mu(z_1+z_3)} \cdot M \cdot e^{-\lambda (z_1+z_3)} dz_1 dz_3 dz_6 dz_8\\ 
    &=& \frac{1}{4} e^{-\mu(2x+2y)} \int_{\mathbb{R}_+^2}  e^{-\mu(z_1+z_3)} \cdot M \cdot e^{-\lambda (z_1+z_3)} dz_1 dz_3\\ 
    &+& \frac{1}{2} \mu e^{-\mu(2x+2y)} \int_{\mathbb{R}_+^3} z_1 \cdot e^{-\mu(z_1+z_3)} \cdot M \cdot e^{-\lambda (z_1+z_3)} dz_1 dz_3 dz_8\\ 
    &=& \frac{1}{4} e^{-\mu(2x+2y)} \left( \frac{\lambda}{\mu + \lambda} \right)^{j+k} \cdot \left( 1+ \frac{ 2\mu k }{\mu + \lambda} \right).
\end{eqnarray*}

\subsection{Condition 2:}
Let $\mathcal{A}_{k}, \mathcal{A}_{j}$ be the inter-arrival times between customer $m$, $m+k$, $n$ and $n+j$ respectively. In addition, we define $M$ to be $$ M = \frac{\lambda^{(k) + (j)} z_1^{k-1} z_2^{j-1}}  {\Gamma(k) \Gamma(j)}.$$This implies that:
\begin{eqnarray*}
    \lefteqn{ \mathbb{P} \left( S^{(1)}_n - \mathcal{A}_{j} > x , S^{(1)}_{n+j} > x, S^{(2)}_{m} > y, S^{(1)}_{m+k} + S^{(2)}_{m+k} + \mathcal{A}_{k} - S^{(1)}_{m} > y , \mathcal{E}_{(2)}^C  \right) } \\  
    &=& \int_{\mathbb{R}_+^8} \{z_4 > x + z_3\} \cdot \{ z_5 > x \} \cdot \{ z_7 > y \} \cdot \{ z_8 + z_9 + z_1 - z_6 > y \} \cdot \{ z_1 + z_8 < z_6 \}\\
    &\times& \mu^6 e^{-\mu(z_4+z_5+z_6+z_7+z_8+z_9)} \cdot M \cdot e^{-\lambda (z_1+z_3)} dz_1 dz_3 dz_4 dz_5 dz_6 dz_7 dz_8 dz_9\\
    &=& e^{-\mu(2x+2y)} \int_{\mathbb{R}_+^4} \{ z_1 + z_8 < z_6 \} \cdot \mu^2 e^{-\mu(-z_1+z_3+2z_6)} \cdot M \cdot e^{-\lambda (z_1+z_3)} dz_1 dz_3 dz_6 dz_8\\
    &=& \frac{1}{4} e^{-\mu(2x+2y)} \int_{\mathbb{R}_+^2}  e^{-\mu(z_1+z_3)} \cdot M \cdot e^{-\lambda (z_1+z_3)} dz_1 dz_3\\
    &=& \frac{1}{4} e^{-\mu(2x+2y)} \left( \frac{\lambda}{\mu + \lambda} \right)^{j+k}.
\end{eqnarray*}

\section{Proof for CASE IX in \ref{ovlp_diff_customers_prop}}
\subsection{Condition 1:}
Let $\mathcal{A}_{n}, \mathcal{A}_{j}, \mathcal{A}_{k}$ be the inter-arrival times between customer $m$, $n$, $n$ and $n+j$, $n+j$ and $m+k$ respectively. In addition, we define $M$ to be $$ M = \frac{\lambda^{(n-m) + (j) + (m+k-n-j)} z_1^{n-m-1} z_2^{j-1} z_3^{m+k-n-j-1}}  {\Gamma(n-m) \Gamma(j) \Gamma(m+k-n-j)}.$$ This implies that:
\begin{eqnarray*}
    \lefteqn{ \mathbb{P} \left( S^{(1)}_n - \mathcal{A}_{j} > x , S^{(1)}_{n+j} > x, S^{(2)}_{m+k} > y, S^{(1)}_{m} + S^{(2)}_{m} - \mathcal{A}_{k} - \mathcal{A}_{n} - \mathcal{A}_{j} - S^{(1)}_{m+k} > y , \mathcal{E}_{(2)}  \right) } \\  
    &=& \int_{\mathbb{R}_+^9} \{z_4 > x + z_2\} \cdot  \{ z_9 > y \} \cdot \{ z_6 + z_7 - z_1 - z_2 - z_3 - z_8 > y \} \\
    &\times& \{ z_1 + z_2 + z_3 + z_8 > z_6 \} \cdot \{ z_5 > x \} \cdot\mu^6 e^{-\mu(z_4+z_5+z_6+z_7+z_8+z_9)} \\
    &\times& M \cdot e^{-\lambda (z_1+z_2+z_3)} dz_1 dz_2 dz_3 dz_4 dz_5 dz_6 dz_7 dz_8 dz_9\\
    &=& e^{-\mu(2x+2y)} \int_{\mathbb{R}_+^5}  \{ z_1 + z_2 + z_3 + z_8 > z_6 \} \\
    &\times& \mu^2 e^{-\mu(z_1+2z_2+z_3+2z_8)} \cdot M \cdot e^{-\lambda (z_1+z_2+z_3)} dz_1 dz_2 dz_3 dz_6 dz_8\\
    &=& \frac{1}{2} e^{-\mu(2x+2y)} \int_{\mathbb{R}_+^4}  \{ z_1 + z_2 + z_3 < z_6 \} \cdot\mu e^{-\mu(z_2+2z_6)} \cdot M \cdot e^{-\lambda (z_1+z_2+z_3)} dz_1 dz_2 dz_3 dz_6\\
    &+& \frac{1}{2} e^{-\mu(2x+2y)} \int_{\mathbb{R}_+^5}  \{ z_1 + z_2 + z_3 > z_6 \} \cdot\mu e^{-\mu(z_1+2z_2+z_3)} \cdot M \cdot e^{-\lambda (z_1+z_2+z_3)} dz_1 dz_2 dz_3 dz_6\\
    &=& \frac{1}{4} e^{-\mu(2x+2y)} \int_{\mathbb{R}_+^3}  e^{-\mu(z_1+2z_2+z_3)} \cdot M \cdot e^{-\lambda (z_1+z_2+z_3)} dz_1 dz_2 dz_3\\
    &+& \frac{1}{2} \mu e^{-\mu(2x+2y)} \int_{\mathbb{R}_+^4}  ( z_1 + z_2 + z_3) \cdot e^{-\mu(z_1+2z_2+z_3)} \cdot M \cdot e^{-\lambda (z_1+z_2+z_3)} dz_1 dz_2 dz_3\\
    &=& \frac{1}{4} e^{-\mu(2x+2y)} \left( \frac{\lambda}{\mu + \lambda} \right)^{k-j} \left( \frac{\lambda}{2\mu + \lambda} \right)^{j} \cdot \left( 1+ 2\mu\left( \frac{k-j}{\lambda+\mu} + \frac{j}{2\mu + \lambda} \right) \right).
\end{eqnarray*}

\subsection{Condition 2:}
Let $\mathcal{A}_{n}, \mathcal{A}_{j}, \mathcal{A}_{k}$ be the inter-arrival times between customer $m$, $n$, $n$ and $n+j$, $n+j$ and $m+k$ respectively. In addition, we define $M$ to be $$ M = \frac{\lambda^{(n-m) + (j) + (m+k-n-j)} z_1^{n-m-1} z_2^{j-1} z_3^{m+k-n-j-1}}  {\Gamma(n-m) \Gamma(j) \Gamma(m+k-n-j)}.$$ This implies that:
\begin{eqnarray*}
    \lefteqn{ \mathbb{P} \left( S^{(1)}_n - \mathcal{A}_{j} > x , S^{(1)}_{n+j} > x, S^{(2)}_{m} > y, -S^{(1)}_{m} + S^{(2)}_{m+k} + \mathcal{A}_{k} + \mathcal{A}_{n} + \mathcal{A}_{j} + S^{(1)}_{m+k} > y , \mathcal{E}_{(2)}^C  \right) } \\  
    &=& \int_{\mathbb{R}_+^9} \{z_4 > x + z_2\} \cdot  \{ z_7 > y \} \cdot \{ -z_6 + z_9 + z_1 + z_2 + z_3 + z_8 > y \} \\
    &\times& \{ z_1 + z_2 + z_3 + z_8 < z_6 \} \cdot \{ z_5 > x \} \cdot\mu^6 e^{-\mu(z_4+z_5+z_6+z_7+z_8+z_9)} \\
    &\times& M \cdot e^{-\lambda (z_1+z_2+z_3)} dz_1 dz_2 dz_3 dz_4 dz_5 dz_6 dz_7 dz_8 dz_9\\
    &=& \frac{1}{4} e^{-\mu(2x+2y)} \int_{\mathbb{R}_+^3}  \mu e^{-\mu(z_1+2z_2+z_3)} \cdot M \cdot e^{-\lambda (z_1+z_2+z_3)} dz_1 dz_2 dz_3\\
    &=& \frac{1}{4} e^{-\mu(2x+2y)} \left( \frac{\lambda}{\mu + \lambda} \right)^{k-j} \left( \frac{\lambda}{2\mu + \lambda} \right)^{j}.
\end{eqnarray*}

\bibliographystyle{plainnat}
\bibliography{tandem}
\end{document}